\documentclass[10pt]{article}
\usepackage{bbm}
\usepackage{amsmath}
\usepackage{amsfonts}
\usepackage{amssymb}
\usepackage{multirow}
\usepackage{mathrsfs}
\usepackage{cite}
\usepackage{ascii}
\usepackage{bbding}
\usepackage{pifont}
\usepackage{ifsym}
\usepackage{wasysym}
\usepackage{upgreek}
\usepackage{color}
\usepackage{graphicx}
\usepackage{eso-pic}
\usepackage{epstopdf}
\usepackage{setspace}
\usepackage{enumerate}
\usepackage{tikz}
\usepackage{algorithm}
\usepackage{algpseudocode}
\usepackage{colortbl}
\usepackage{booktabs}
\usepackage{makecell}
\usepackage{appendix}
\usepackage{array}
\usepackage{threeparttable}
\usepackage{diagbox}

\renewcommand{\Box}{\framebox{\rule{0.3em}{0.0em}}}
\newtheorem{thm}{Theorem}[section]

\newtheorem{ex}{Example}[section]


\renewcommand{\Box}{\hfill\rule{2.3mm}{2.3mm}}

\renewcommand{\Box}{\framebox{}}
\numberwithin{equation}{section}
\title{Theoretical and numerical comparison of seven single-level reformulations for bilevel programs 
}
\author{Yu-Wei Li\thanks{\baselineskip 9pt
School of Management, Shanghai University, Shanghai 200444, China. E-mail: yuwei\_li@shu.edu.cn.}, \
   Gui-Hua Lin\thanks{\baselineskip 9pt
Corresponding author. School of Management, Shanghai University, Shanghai 200444, China. E-mail:
guihualin@shu.edu.cn.}, \
        Xide Zhu\thanks{\baselineskip 9pt
School of Management, Shanghai University, Shanghai 200444, China. E-mail: xidezhu@shu.edu.cn.}
}

\date{February 9, 2025}

\begin{document}
\maketitle

\baselineskip 16pt

\vspace{4pt} \noindent{\bf Abstract.}
This paper considers a bilevel program. To solve this bilevel program, it is generally necessary to transform it into some single-level optimization problem. One approach is to replace the lower-level program by its KKT conditions to transform the bilevel program as a mathematical program with complementarity constraints (MPCC). Another approach is to apply the lower-level Wolfe/Mond-Weir/extended Mond-Weir duality to transform the bilevel program into some duality-based single-level reformulations, called WDP, MDP, and eMDP respectively in the literature. In this paper, inspired by a conjecture from a recent publication that the tighter feasible region of a reformulation, the better its numerical performance, we present three new duality-based single-level reformulations, called TWDP/TMDP/eTMDP, with tighter feasible regions. Our main goal is to compare all above-mentioned reformulations by designing some direct and relaxation algorithms with projection and implementing these algorithms on 450 test examples generated randomly. Our numerical experiments show that, whether overall comparison or pairwise comparison, at least in our tests, the WDP/MDP/TWDP/TMDP reformulations were always better than the MPCC reformulation, while the eMDP/eTMDP reformulations were always the worst ones among six duality-based reformulations, which indicates that the above conjecture is incorrect. In particular, for the relaxation algorithms, the WDP/MDP/TWDP/TMDP reformulations performed 3-5 times better than the MPCC reformulation, while the eMDP/eTMDP reformulations performed 2 times better than the MPCC reformulation.

\vspace{4pt}\noindent{\bf Keywords.} Bilevel program, Wolfe duality, Mond-Weir duality, MPCC, direct/relaxation algorithms.

\vspace{4pt}\noindent{\bf 2010 Mathematics Subject Classification.} \ 90C30, 90C33, 90C46
\baselineskip 18pt
\parskip 2pt

\section{Introduction}

In this paper, we consider the bilevel program
\begin{eqnarray*}\label{BP}
{\rm( BP)}\qquad\min\limits_{x,y} &&F(x, y) \nonumber\\
\mbox{s.t.}&&(x,y)\in \Omega,~y \in {S}(x),
\end{eqnarray*}
where $x\in\mathbb{R}^{n}$ and $y\in\mathbb{R}^{m}$ are respectively called upper-level and lower-level variables 
and ${S}(x) \subseteq \mathbb{R}^m$ denotes the optimal solution set of the lower-level parameterized program
\begin{eqnarray*}\label{P_x}
{(\mathrm{P}_x)}\qquad\min\limits_{y}&&f(x, y) \nonumber\\
\mbox{s.t.}&&g(x, y) \le 0, ~h(x, y)=0,
\end{eqnarray*}
where $g$: $\mathbb{R}^{n+m} \rightarrow\mathbb{R}^{p}$ and $h$: $\mathbb{R}^{n+m} \rightarrow\mathbb{R}^{q}$. Throughout the paper, we denote by $$X=\{x\in \mathbb{R}^n: \exists~ y\in \mathbb{R}^m~{\rm s.t.}~(x,y)\in\Omega\}, \quad Y(x) = \{y \in \mathbb{R}^{m} : g(x, y) \le 0, ~h(x, y) = 0\}.$$ In addition, we assume that $F$ is continuously differentiable, $\{f, g, h\}$ are all twice continuously differentiable, and ${S}(x) \neq \emptyset$ for each $x\in X$ so as to ensure the bilevel program is always well-defined.

The above bilevel program is not only widely used in classical fields such as principal-agent problems, transportation planning, and supply chain management \cite{Colson2007overview, Khorramfar2022managing, Goyal2023decision}, but also plays a significant role in recent popular fields such as hyperparameter optimization and meta-learning \cite{Liu2023hierarchical, Ye2023difference, Zhang2024introduction}. See, e.g., \cite{Colson2007overview, Sabach2017first, Dempe2020bilevel, Zeng2020practical, Kleinert2021computing, Hong2023two, Lu2024first} and references therein for more details.

Due to its hierarchical structure, the bilevel program BP is very difficult to solve. Actually, it is NP-hard even if all functions involved are linear \cite{Ben1990computational}. Understandably, to deal with this problem, it generally needs to transform it into some single-level optimization problems so that one can apply the well-developed optimization algorithms in mathematical programming. Along this way,
the most popular approach is to utilize the Karush-Kuhn-Tucker (KKT) conditions of the lower-level program to transform BP into the following mathematical program with complementarity constraints:
\begin{eqnarray*}\label{MPCC}
{\rm(MPCC)}\qquad\min && F(x,y) \\
\mbox{s.t.} && (x,y)\in \Omega,~\nabla_y f(x, y) + \nabla_y g(x, y) u+ \nabla_y h(x, y)v=0,\\
&&u\geq0,~g(x, y)\le 0,~ u^Tg(x, y)=0, ~h(x, y)=0.
\end{eqnarray*}
However, one well-known flaw of the above approach is that MPCC fails to satisfy the Mangasarian-Fromovitz constraint qualification (MFCQ) at any feasible point so that the well-developed optimization algorithms may not be stable in solving MPCC.
After decades of development, numerous approximate methods such as relaxation methods, penalty function methods, interior point methods, implicit programming methods, active-set identification methods, constrained equation-based methods, etc., have been proposed for solving MPCC.
See, e.g., \cite{Luo1996mathematical, Leyffer2006interior, Lin2006hybrid, Lin2009solving, Izmailov2012semismooth, LinEquiation,Scholtes2001} and references therein for more details above the developments on MPCC.

Recently, we present a new approach to deal with BP based on the lower-level Wolfe or Mond-Weir duality \cite{Li2023novel, Li2024solving}. This approach contributes the following three new single-level reformulations, which are based on the lower-level Wolfe/Mond-Weir/extended Mond-Weir duality respectively:
\begin{eqnarray*}\label{WDP}
{\rm (WDP)}\qquad\min &&F(x, y)  \\
\mbox{s.t.} && (x,y)\in \Omega,~g(x, y) \le 0, ~h(x, y)=0,\\
&&f(x, y) - f(x, z) - u^T g(x, z) - v^T h(x, z) \leq 0, \\
&&\nabla_z f(x, z) + \nabla_z g(x, z) u + \nabla_z h(x, z) v = 0, ~u\geq 0;
\end{eqnarray*}
\begin{eqnarray*}\label{MDP}
{\rm (MDP)}\qquad\min &&F(x, y) \nonumber \\
\mbox{\rm s.t.} &&(x,y)\in \Omega,~g(x, y) \le 0, ~h(x, y)=0,\\
&&f(x, y) - f(x, z)\leq 0,~ u^T g(x, z) + v^T h(x,z)\geq0,\nonumber\\
&&\nabla_z f(x, z) + \nabla_z g(x, z) u + \nabla_z h(x, z) v = 0, ~u\geq 0;\nonumber
\end{eqnarray*}
\begin{eqnarray*}\label{eMDP}
~~~~{\rm (eMDP)}\qquad\min &&F(x, y) \nonumber \\
\mbox{\rm s.t.} &&(x,y)\in \Omega,~g(x, y) \le 0, ~h(x, y)=0,\\
&&f(x, y) - f(x, z)\leq 0,~ u\circ g(x, z) \geq0, ~v \circ h(x,z)=0,\nonumber\\
&&\nabla_z f(x, z) + \nabla_z g(x, z) u + \nabla_z h(x, z) v = 0, ~u\geq 0.\nonumber
\end{eqnarray*}
Here, $\circ$ denotes the Hadamard product. They have two advantages over MPCC: (i) In theory, WDP/MDP/eMDP may satisfy the MFCQ at their feasible points, which reveals that the new ones have theoretical and structural advantages over MPCC. (ii) Numerous numerical examples show that WDP/MDP/eMDP are more effective than MPCC in solving the original bilevel programs. See \cite{Li2023novel, Li2024solving,Li-PJO} for more details. It needs to point out that another remarkable approach to deal with BP is based on the lower-level optimal value function \cite{Ye2010new, Mehlitz2021sufficient,Lin2014simple, Xu2014smoothing}. Since the optimal value function does not have analytic expression generally, it is hard to compare this approach with the ones without extremal functions introduced above.

In this paper, we continue to study the approach based on lower-level duality and, particularly, we will propose three new single-level reformulations for BP. Our motivation comes from an observation from our work \cite{Li2024solving}. That is, through 150 numerical examples generated randomly in \cite{Li2024solving}, it seems that the tighter the feasible region of the reformulation, the better the corresponding numerical performance. Inspired by this observation, we may move the equality constraint $h(x,z)=0$ outside from WDP/MDP/eMDP so that the feasible regions become tighter than before. Thus, we can obtain three new reformulations for BP, which are called TWDP/TMDP/eTMDP below respectively. Naturally, in order to ensure the equivalence between the new reformulations and the original BP, some duality theorems need to be established firstly. See Section \ref{sec-NewReformulations} for details.

Main contributions of this paper can be summarized as follows:
\begin{itemize}
  \item Three new single-level reformulations TWDP/TMDP/eTMDP are proposed for BP and new duality theorems are established. Based on these duality theorems, TWDP/TMDP/eTMDP are shown to be equivalent to BP under mild conditions in the globally or locally optimal sense. Thus, there are seven single-level reformulations for BP so far, including the old ones MPCC/WDP/MDP/eMDP and the new ones TWDP/TMDP/eTMDP.

  \item The most important goal in this paper is to present overall and pairwise comparison of the above-mentioned seven single-level reformulations by numerical experiments. Since a lot of numerical examples in \cite{Li2023novel,Li2024solving} can not be feasibly solved by the proposed algorithms, as an improvement, we add some projection steps into the algorithms of this paper so as to ensure to output feasible points of BP. We further establish comprehensive convergence analysis for these algorithms.

  \item We randomly generated 450 test examples of three types of bilevel programs, which include 150 bilevel programs with LP constraints, 150 bilevel programs with QP constraints, and 150 bilevel programs with QCQP constraints. Our numerical experiments show that, whether overall comparison or pairwise comparison, at least in our tests, the WDP/MDP/TWDP/TMDP reformulations were always better than the MPCC reformulation, while the extended two reformulations eMDP and eTMDP were always the worst ones among the duality-based reformulations, which shows that the above-mentioned conjecture from \cite{Li2024solving} is incorrect. In particular, the WDP/MDP/TWDP/TMDP reformulations were 3-5 times better than the MPCC reformulation, while the eMDP/eTMDP reformulations were 2 times better than the MPCC reformulation.
\end{itemize}

The rest of the paper is organized as follows. In Section \ref{sec-NewReformulations}, we introduce the new single-level reformulations TWDP/TMDP/eTMDP, give some useful duality theorems, show their equivalence with BP in the globally and locally optimal senses, and provide some examples to show that the new reformulations may satisfy the MFCQ at their feasible points. In Section 3, we investigate the relations among seven single-level reformulations. In Section \ref{sec-algorithms}, we present
two algorithms with projection, called direct and relaxation algorithms respectively. We also give some convergence results for these algorithms. In Section \ref{sec-Numerical}, we give numerical comparison by implementing the proposed algorithms on 450 test examples generated randomly, which includes a comparison of direct and relaxation algorithms, an overall comparison of seven single-level reformulations, and a pairwise comparison of seven single-level reformulations. From these numerical comparisons, we can obtain their own advantages and disadvantages. Some concluding remarks are given in Section \ref{sec-conclusion}. To improve readability, all proofs are given in the online supplement as appendices and, to avoid excessive length, all numerical results for 450 test examples are provided in the data repository.

\section{New reformulations with tighter feasible regions }
\label{sec-NewReformulations}

As stated in Section 1, motivated by the conjecture from \cite{Li2024solving} that the reformulation with tighter feasible region may have better numerical performance, we may move the equality constraint $h(x,z)=0$ outside from WDP/MDP/eMDP so that the feasible regions become tighter than before. Then, we can obtain the following new single-level reformulations for BP:
\begin{eqnarray*}\label{TWDP}
{\rm (TWDP)}\qquad\min &&F(x, y)  \\
\mbox{s.t.} && (x,y)\in \Omega,~g(x, y) \le 0, ~h(x, y)=0,\\
&&f(x, y) - f(x, z) - u^T g(x, z) \leq 0,~h(x, z) =0, \\
&&\nabla_z f(x, z) + \nabla_z g(x, z) u + \nabla_z h(x, z) v = 0, ~u\geq 0;
\end{eqnarray*}
\begin{eqnarray*}\label{TMDP}
{\rm (TMDP)}\qquad\min &&F(x, y) \nonumber \\
\mbox{\rm s.t.} &&(x,y)\in \Omega,~g(x, y) \le 0, ~h(x, y)=0,\\
&&f(x, y) - f(x, z)\leq 0,~ u^T g(x, z) \geq0,~h(x, z) =0,\nonumber\\
&&\nabla_z f(x, z) + \nabla_z g(x, z) u + \nabla_z h(x, z) v = 0, ~u\geq 0;\nonumber
\end{eqnarray*}
\begin{eqnarray*}\label{eTMDP}
{\rm (eTMDP)}\qquad\min &&F(x, y) \nonumber \\
\mbox{\rm s.t.} &&(x,y)\in \Omega,~g(x, y) \le 0, ~h(x, y)=0,\\
&&f(x, y) - f(x, z)\leq 0,~ u\circ g(x, z) \geq0, ~h(x, z) =0,\nonumber\\
&&\nabla_z f(x, z) + \nabla_z g(x, z) u + \nabla_z h(x, z) v = 0, ~u\geq 0.\nonumber
\end{eqnarray*}
In order to ensure the equivalence between the new reformulations and BP, we need to establish some corresponding duality theorems.

\subsection{Duality theorems}

Responding to the new reformulations given above, we present the following three new dual problems for the lower-level program $\mathrm{P}_x$:
\begin{eqnarray*}
(\mathrm{TWD}_x)\qquad\max\limits_{z, u, v}&&f(x, z) + u^T g(x, z) \nonumber\\
\mbox{s.t.}&& \nabla_z L(x, z, u, v)=0,~h(x,z)=0,~u\geq 0;\qquad\qquad\qquad
\end{eqnarray*}
\begin{eqnarray*}
(\mathrm{TMD}_x)\qquad\max\limits_{z, u, v}&&f(x, z) \nonumber\\
\mbox{s.t.}&& \nabla_z L(x, z, u, v)=0,~u^T g(x, z)\geq0,~h(x,z)=0,~u\geq 0;
\end{eqnarray*}
\begin{eqnarray*}
(\mathrm{eTMD}_x)\qquad\max\limits_{z, u, v}&&f(x, z) \nonumber\\
\mbox{s.t.}&& \nabla_z L(x, z, u, v)=0,~u\circ g(x, z)\geq0,~h(x,z)=0,~u\geq 0.
\end{eqnarray*}
Here, $L(x, z, u, v) = f(x, z) + u^T g(x, z) + v^T h(x, z).$ These dual problems are generalizations of the Wolfe duality in \cite{Wolfe1961duality}, the Mond-Weir duality in \cite{Mond1981generalized}, and the extended Mond-Weir duality in \cite{Li2024solving} respectively.
Particularly, if $\mathrm{P}_x$ does not contain equality constraints, these new ones reduce to the original Wolfe/Mond-Weir/extended Mond-Weir dual problems.

Denote by $\mathcal{F}_1(x), \mathcal{F}_2(x), \mathcal{F}_3(x)$ the feasible regions of $\mathrm{TWD}_x$, $\mathrm{TMD}_x$, $\mathrm{eTMD}_x$ respectively and let
\begin{eqnarray*}
\mathcal{Z}_1(x,u,v)&=&\{z\in \mathbb{R}^m: ~\nabla_z L(x,z, u, v)=0,~g(x,z)\leq0,~h(x,z)=0\},\\
\mathcal{Z}_2(x,u,v)&=&\{z\in \mathbb{R}^m: ~\nabla_z L(x,z, u, v)=0,~u^T g(x, z)\geq0,~g(x,z)\leq0,~h(x,z)=0\},\\
\mathcal{Z}_3(x,u,v)&=&\{z\in \mathbb{R}^m: ~\nabla_z L(x,z, u, v)=0,~u\circ g(x, z)\geq0,~g(x,z)\leq0,~h(x,z)=0\}.
\end{eqnarray*}
Then, we have the following weak and strong duality theorems, which are shown in Appendix A.

\begin{thm}\label{TWD-dualitythm}\emph{\textbf{(Wolfe case)}}
Assume that, for any given $x\in X$,  $u\in\mathbb{R}^p_{+}$, and $v\in\mathbb{R}^q$, $L(x, \cdot, u, v)$ is pseudoconvex on $\mathcal{Z}_1(x,u,v)$. Then, the weak duality holds, i.e.,
\begin{eqnarray*}\label{mWD-dualitythm-1}
\min\limits_{y\in Y(x)} f(x,y)  \geq \max\limits_{(z, u, v)\in \mathcal{F}_1(x)} f(x,z)+u^Tg(x,z).
\end{eqnarray*}
Moreover, if $\mathrm{P}_x$ satisfies the Guignard constraint qualification {\rm (GCQ)} at some solution $y_x\in S(x)$, there exists $(z_x,u_x,v_x)\in \mathcal{F}_1(x)$ such that the strong duality holds, i.e.,
\begin{eqnarray*}\label{mWD-dualitythm-2}
\min\limits_{y\in Y(x)} f(x,y)  =f(x,y_x)=f(x,z_x)+u_x^Tg(x,z_x)= \max\limits_{(z, u, v)\in \mathcal{F}_1(x)} f(x,z)+u^Tg(x,z).
\end{eqnarray*}
\end{thm}

\begin{thm}\label{TMD-dualitythm}\emph{\textbf{(Mond-Weir case)}}
Assume that, for any given $x\in X$, $u\in\mathbb{R}^p_{+}$, and $v\in\mathbb{R}^q$, $f(x,\cdot)$ is pseudoconvex and $u^Tg(x,\cdot)+v^Th(x,\cdot)$ is quasiconvex on $\mathcal{Z}_2(x,u,v)$. Then, the weak duality holds, i.e.,
\begin{eqnarray*}\label{mMD-dualitythm-1}
\min\limits_{y\in Y(x)} f(x,y)  \geq \max\limits_{(z, u, v)\in \mathcal{F}_2(x)} f(x,z).
\end{eqnarray*}
Moreover, if $\mathrm{P}_x$ satisfies the {\rm GCQ} at some solution $y_x\in S(x)$,  the strong duality holds, i.e.,
\begin{eqnarray*}\label{mMD-dualitythm-2}
\min\limits_{y\in Y(x)} f(x,y)  =f(x,y_x)= \max\limits_{(z, u, v)\in \mathcal{F}_2(x)} f(x,z).
\end{eqnarray*}
\end{thm}

\begin{thm}\label{meMD-dualitythm}\emph{\textbf{(extended Mond-Weir case)}}
Assume that, for any given $x\in X$, $u\in\mathbb{R}^p_{+}$, and $v\in\mathbb{R}^q$, $f(x,\cdot)$ is pseudoconvex and $u^Tg(x,\cdot)+v^Th(x,\cdot)$ is quasiconvex on $\mathcal{Z}_3(x,u,v)$. Then, the weak duality holds, i.e.,
\begin{eqnarray*}\label{meMD-dualitythm-1}
\min\limits_{y\in Y(x)} f(x,y)  \geq \max\limits_{(z, u, v)\in \mathcal{F}_3(x)} f(x,z).
\end{eqnarray*}
Moreover, if $\mathrm{P}_x$ satisfies the {\rm GCQ} at some solution $y_x\in S(x)$,  the strong duality holds, i.e.,
\begin{eqnarray*}\label{me9MD-dualitythm-2}
\min\limits_{y\in Y(x)} f(x,y)  =f(x,y_x)= \max\limits_{(z, u, v)\in \mathcal{F}_3(x)} f(x,z).
\end{eqnarray*}
\end{thm}

\subsection{Equivalence between TWDP/TMDP/eTMDP and BP}

With the help of the above weak and strong duality theorems, we can obtain the following equivalence results between TWDP/TMDP/eTMDP and BP.

\begin{thm}
Suppose that, for any $x\in X$, $u\in \mathbb{R}^p_+$, and $v\in\mathbb{R}^q$, $L(x, \cdot, u, v)$ is pseudoconvex on $\mathcal{Z}_1(x,u,v)$ in the case of {\rm TWDP}, or $f(x,\cdot)$ is pseudoconvex and $u^T g(x,\cdot)+v^Th(x,\cdot)$ is quasiconvex on $\mathcal{Z}_i(x,u,v)~(i=2,3)$ in the cases of {\rm TMDP/eTMDP}. Suppose also that $\mathrm{P}_x$ satisfies the GCQ at some solution $y_x \in {S(x)}$.
Then, if $(\bar{x},\bar{y})$ is an optimal solution to {\rm BP}, there exists $(\bar{z},\bar{u},\bar{v})\in \mathcal{F}_1(\bar{x})$~$(or~\mathcal{F}_2(\bar{x})/\mathcal{F}_3(\bar{x}))$ such that $(\bar{x},\bar{y},\bar{z},\bar{u},\bar{v})$ is an optimal solution to {\rm TWDP~(}or {\rm TMDP/eTMDP)}.
Conversely, if $(\bar{x},\bar{y},\bar{z},\bar{u},\bar{v})$ is an optimal solution to {\rm TWDP~(}or {\rm TMDP/eTMDP)}, then $(\bar{x},\bar{y})$ is an optimal solution to {\rm BP}.
\end{thm}

\begin{thm}
If $\mathrm{P}_x$ and $\mathrm{TWD}_x$ {\rm(}or $\mathrm{TMD}_x/\mathrm{eTMD}_x${\rm)} satisfy the weak duality and $(x, y, z, u, v)$ is feasible to {\rm TWDP~(}or {\rm TMDP/eTMDP)}, then $y$ is globally optimal to $\mathrm{P}_x$ and $(z, u, v)$ is globally optimal to $\mathrm{TWD}_x$ {\rm(}or $\mathrm{TMD}_x/\mathrm{eTMD}_x${\rm)}.
\end{thm}

\begin{thm}
Let $(\bar{x},\bar{y},\bar{y},\bar{u},\bar{v}) $ be locally optimal to {\rm TWDP~(}or {\rm TMDP/eTMDP)} with $\{\bar{u},\bar{v}\}$ being the corresponding Lagrange multipliers. Under any conditions to ensure the weak duality for $\mathrm{P}_{\bar{x}}$ and $\mathrm{TWD}_{\bar{x}}$ {\rm(}or $\mathrm{TMD}_{\bar{x}}/\mathrm{eTMD}_{\bar{x}}${\rm)},
if $\mathrm{P}_x$ satisfies the MFCQ at $y\in S(x)$ for all $(x,y)$ closely to $(\bar{x},\bar{y})$, $(\bar{x},\bar{y})$ is locally optimal to {\rm BP}.
\end{thm}

Since the proofs of these theorems are similar to the corresponding results given in \cite{Li2023novel,Li2024solving}, we omit them here.

\subsection{TWDP/TMDP/eTMDP may satisfy the MFCQ at feasible points}

It is well-known that MPCC fails to satisfy the MFCQ at any feasible point. Unlike MPCC, it has been shown in \cite{Li2023novel,Li2024solving} that the lower-level duality-based reformulations WDP/MDP/eMDP may satisfy the MFCQ at their feasible points. The following examples indicate that the new reformulations TWDP/TMDP/eTMDP may also satisfy the MFCQ at their feasible points.

\begin{ex}\rm\label{ex-mWDP-MFCQ}
Consider the bilevel program
\begin{eqnarray}\label{ex-mWDP-MFCQ-1}
\min &&(x-y_1-8)^2\\
\mbox{s.t.}&&x \geq 1 ,~y \in {S}(x)=\arg\min\limits_{y}\{y_1-y_2:y_1^3\leq x,~ y_1\geq0,~y_1+y_2-x=0\}.\nonumber
\end{eqnarray}
It is easy to see that $S(x)=\{(0,x)\}$ for any $x \geq 1$ and then \eqref{ex-mWDP-MFCQ-1} has a unique minimizer $(8,0,8)$. It is not difficult to verify that $(x^*, y^*, z^*, u^*,v^*)=(8, 0, 8,-2,10,0, 2, 1)$ is a feasible and global optimal solution to its TWDP reformulation
\begin{eqnarray}\label{ex-mWDP-MFCQ-2}
\min&&(x-y_1-8)^2 \nonumber\\
\mbox{s.t.}&&y_1-y_2-z_1+z_2-u_1(z_1^3-x)+u_2z_1\leq0, \nonumber\\
&&1+3z_1^2u_1-u_2+v=0,~-1+v=0,\\
&&u_1\geq0,~u_2\geq0,~x \geq 1,~y_1^3\leq x,~ y_1\geq0,\nonumber\\
&&y_1+y_2-x=0,~z_1+z_2-x=0.\nonumber
\end{eqnarray}
By direct calculation, we can find a vector $d=(1,1,0,0,1,1,12,0)^T$ such that $d^T\mathbf{h}_i=0~(i=1,\cdots,4)$ and $d^T\mathbf{g}_i<0~(i=1,2,3)$, where $\mathbf{h}_i~(i=1,\cdots,4)$ denote the gradients corresponding to equality constraints and $\mathbf{g}_i~(i=1,2,3)$ denote the gradients corresponding to active inequality constraints at $(8, 0, 8,-2,10,0, 2, 1)$. This means that \eqref{ex-mWDP-MFCQ-2} satisfies the MFCQ at $(8, 0, 8,-2,10,0, 2, 1)$.
%
%
\end{ex}

\begin{ex}\rm\label{ex-mMDP-MFCQ}
Consider the bilevel program
\begin{eqnarray}\label{ex-mMDP-MFCQ-1}
\min &&(x+y_1+y_2)^2\\
\mbox{s.t.}&&x \in[-1,1] ,~y \in {S}(x)=\arg\min\limits_{y}\{y_1^3-3y_1:y_1\geq x,~y_1+y_2-x=0\}.\nonumber
\end{eqnarray}
Obviously, $S(x)=\{(1,x-1)\}$ for any $x\in[-1,1]$ and \eqref{ex-mMDP-MFCQ-1} has a unique minimizer $(0,1,-1)$.
The TMDP/eTMDP reformulations for \eqref{ex-mMDP-MFCQ-1} are the same problem
\begin{eqnarray}\label{ex-mMDP-MFCQ-2}
\min&&(x+y_1+y_2)^2 \nonumber\\
\mbox{s.t.}&&y_1^3-3y_1-z_1^3+3z_1\leq0,~u(x-z_1)\geq0, \nonumber\\
&&3z_1^2-3-u+v=0,~v=0,\\
&&u\geq0,~x \in[-1,1],~y_1\geq x,\nonumber\\
&&y_1+y_2-x=0,~z_1+z_2-x=0.\nonumber
\end{eqnarray}
It is easy to verify that $(x^*, y^*, z^*, u^*,v^*)=(0, 1, -1, -2, 2, 9, 0)$ is a feasible and global optimal solution to \eqref{ex-mMDP-MFCQ-2}. Let $\mathbf{h}_i~(i=1,\cdots,4)$ be the gradients corresponding to equality constraints and $\mathbf{g}$ be the gradient corresponding to the unique active inequality constraint at this point. Then, we can find a vector $d=(0,0,0,1,-1,-12,0)^T$ such that $d^T\mathbf{h}_j=0~(j=1,\cdots,4)$ and $d^T\mathbf{g}<0$, which means that \eqref{ex-mMDP-MFCQ-2} satisfies the MFCQ at $(0, 1, -1, -2, 2, 9, 0)$.
%
%
\end{ex}

\section{Relationship among seven single-level reformulations}
\label{sec-Comparison existing single-level reformulations}
So far, there have been proposed seven single-level reformulations with analytical expressions for BP, which include the old ones MPCC/WDP/MDP/eMDP and the new ones TWDP/TMDP/eTMDP. The following result related to these reformulations is evident.

%
%


\begin{thm}\label{Relations-feasibility}
The point $(x, y, u, v)$ is feasible to {\rm MPCC} if and only if $(x, y, y, u,v)$ is feasible to any one of {\rm WDP/MDP/eMDP/TWDP/TMDP/eTMDP}. Therefore, the duality-based reformulations can be regarded as dimensionality extensions from MPCC on lower-level variables.
Moreover, their feasible regions have the following relationship:
\begin{enumerate}[\rm (i)]
\item $\mathcal{F}_{\rm \scriptscriptstyle eMDP}\subseteq\mathcal{F}_{\rm \scriptscriptstyle MDP}\subseteq\mathcal{F}_{\rm \scriptscriptstyle WDP}$, \ $\mathcal{F}_{\rm \scriptscriptstyle eTMDP}\subseteq\mathcal{F}_{\rm \scriptscriptstyle TMDP}\subseteq\mathcal{F}_{\rm \scriptscriptstyle TWDP}$;
\item $\mathcal{F}_{\rm \scriptscriptstyle TWDP}\subseteq\mathcal{F}_{\rm \scriptscriptstyle WDP}$, $\mathcal{F}_{\rm \scriptscriptstyle TMDP}\subseteq\mathcal{F}_{\rm \scriptscriptstyle MDP}$, $\mathcal{F}_{\rm \scriptscriptstyle eTMDP}\subseteq\mathcal{F}_{\rm \scriptscriptstyle eMDP}$.
\end{enumerate}
In particular, among the duality-based reformulations, the feasible region of WDP is the largest, while the feasible region of eTMDP is the tightest.
\end{thm}

On the other hand, since all above-mentioned single-level reformulations are nonconvex optimization problems, it is necessary to study the relationship among their stationary points. We have the following results, whose proofs are provided in  Appendix B.

\begin{thm}\label{StationaryPoint}
\begin{enumerate}[\rm (i)]
  \item If the point $(\bar{x},\bar{y},\bar{z},\bar{u},\bar{v})$ is a KKT point of {\rm WDP} and  feasible to any one of {\rm MDP/eMDP/TWDP/TMDP/eTMDP}, it is also a KKT point of this problem.
  \item If the point $(\bar{x},\bar{y},\bar{z},\bar{u},\bar{v})$ is a KKT point of {\rm MDP} and  feasible to any one of {\rm eMDP/TMDP/eTMDP}, it is also a KKT point of this problem.
  \item If the point $(\bar{x},\bar{y},\bar{z},\bar{u},\bar{v})$ is a KKT point of {\rm TWDP} and  feasible to any one of {\rm eMDP/TMDP/eTMDP}, it is also a KKT point of this problem.
\end{enumerate}
\end{thm}

From Theorem \ref{StationaryPoint}, it seems that the larger the reformulation's feasible region, the stronger its KKT conditions. However, by the proof of Theorem \ref{StationaryPoint}, we can find that the KKT points of eMDP/TMDP/eTMDP are actually equivalent to each other under feasibility conditions.

\begin{thm}\label{KKTandS}
If the point $(\bar{x},\bar{y},\bar{y},\bar{u},\bar{v})$ is a KKT point of any one of {\rm WDP/MDP/eMDP/TWDP/TMDP/eTMDP}, then $(\bar{x}, \bar{y}, \bar{u},\bar{v})$ is a strongly stationary point of {\rm MPCC}.
\end{thm}

Note that the converse of Theorem \ref{KKTandS} is not true, which can be illustrated by the examples
\begin{eqnarray*}\label{ex-eMDP-KKT-S-1}
\min&& x^2-(2y+1)^2\\
\mbox{s.t.}&& x\leq0,~y\in \arg\min\limits_{y}\{(y-1)^2: 3x-y-3\leq0,~x+y-1\leq0\}\nonumber
\end{eqnarray*}
and
\begin{eqnarray*}\label{ex-TWDPTMDPeTMDP-kkt-s-1}
\min &&x^2-(2y_1+1)^2 \nonumber\\
\mbox{s.t.}&& x\leq0,~y\in  \arg\min\limits_{y}\Big\{\frac{(y_1-1)^2}{2}-x^3y_2: x+\frac{y_1}{2}-y_2-4\leq0,~1-y_2\leq0,~y_1+2y_2+x-3=0\Big\}.
\end{eqnarray*}

\section{Direct and relaxation algorithms with projection}\label{sec-algorithms}
In this section, we give a comprehensive comparison of seven single-level reformulations with analytical expressions for BP through numerical experiments. As in \cite{Li2023novel, Li2024solving}, we focus on comparing their numerical performance in solving the original BP by designing two kinds of algorithms: One uses the direct strategy and the other uses some relaxation strategy. Different from the algorithms in \cite{Li2023novel, Li2024solving}, we add some projection steps into the algorithms as an improvement so as to ensure to output feasible points of BP always in theory at least.

To simplify our discussion, we refer to all single-level reformulations as problem Q and the corresponding relaxation problems as Q$(t)$ with a relaxation parameter $t$. We also simplify the upper-level constraint $(x,y)\in\Omega$ as $x\in X$ in the remaining sections, where $X$ is a nonempty closed convex set. In addition, we uniformly represent the variable of problem Q as $(x, y, z, u, v)$ and, as an exception, the variable in the case of MPCC is actually $(x, y, u, v)$. To ensure the projection onto the lower-level program P$_x$ to be well-defined, we assume that P$_x$ is always convex for any $x\in X$.

How to compare these different single-level reformulations for BP is an important and crucial issue. The fairest way is to evaluate their numerical performance by the original BP. Specifically, we evaluate different reformulations in three levels: The top priority is to check the approximation solutions' feasibilities to the original BP, the second priority is to compare their upper-level objective values if their feasibilities satisfy given accuracy tolerances, and the third priority is to compare CPU times when the top and second priorities are the same.
To this end, we adopt the following consistent criterion to measure the infeasibility to the original BP:
\begin{eqnarray}\label{Inf}
{\tt Infeasibility}=\min_{x\in X}\|x^k-x\|+\|\max(0,g(x^k,y^k))\|+\|h(x^k,y^k)\|+|f(x^k,y^k)-V(x^k)|,
\end{eqnarray}
where $(x^k,y^k)$ is an approximation solution and $V(x^k)=\min \{f(x^k,y): g(x^k,y)\leq 0,~h(x^k,y)=0\}$. Since $X$ is assumed to be a nonempty closed convex set and the lower-level problem P$_x$ is assumed to be convex for any $x\in X$ , the measure `\rm{Infeasibility}' is well-defined.
Obviously, {\tt Infeasibility} equals to zero if and only if $(x^k,y^k)$ is feasible to the original BP.

\subsection{Direct algorithm with projection}

Direct algorithm means to solve problem Q directly to get a candidate solution. If this point is feasible to the original BP, we output it as an approximation solution. Otherwise, we make a projection onto the feasible region of BP. The procedures are described in Algorithm \ref{alg:One Projection Direct Algorithm}.

\begin{algorithm}
  \caption{\quad Direct Algorithm with Projection}
  \label{alg:One Projection Direct Algorithm}

\textbf{Step 0}  Choose an initial point $x^0 \in X$ and two accuracy tolerances $\{\epsilon_{\rm \scriptscriptstyle SQP},\epsilon_{\rm \scriptscriptstyle Inf}\}$.

\textbf{Step 1} 
Solve the lower-level program P$_{x^0}$ to get $({y}^0,{u}^0,{v}^0)$. Set ${w}^0=({x}^0,{y}^0,{y}^0,{u}^0,{v}^0)$.

\textbf{Step 2} Solve problem Q by \emph{fmincon} with ${w}^0$ as a starting point and $\epsilon_{\rm \scriptscriptstyle SQP}$ as an accuracy tolerance to obtain a candidate $\tilde{w}=(\tilde{x},\tilde{y},\tilde{z},\tilde{u},\tilde{v})$.

\textbf{Step 3}
Calculate {\tt Infeasibility} with $\tilde{w}$.
If {\tt Infeasibility} $\leq \epsilon_{\rm \scriptscriptstyle Inf}$, stop. Otherwise, solve $x^*=\arg\min_{x\in X}\|x-\tilde{x}\|^2$ and solve the lower-level program P$_{x^*}$ to obtain a new candidate $(x^*,y^*,y^*,u^*,v^*)$.
\end{algorithm}

When $X$ is  nonempty closed convex and the lower-level programs are convex,
Algorithm \ref{alg:One Projection Direct Algorithm}
must output a feasible solution to the original BP as long as the lower-level programs are solvable.
In the case where some lower-level program has no solution, we need to change an initial point to try again.

\subsection{Relaxation algorithm with projection}

All publications on MPCC and our works \cite{Li2023novel, Li2024solving} about the duality-based approach reveal that solving the single-level reformulation Q directly is not an effective method in the study of bilevel programs. For this sake, some approximation techniques are usually necessary. One popular technique is to approximate problem Q by relaxing those excessive constraints appropriately. For MPCC, we generally relax the constraint $u^T g(x,y)=0$ (which is equivalent to $u^T g(x,y)\ge 0$) by $u^T g(x,y)\ge -t$, where $t>0$ is a small relaxation parameter. For the  duality-based reformulations, we generally relax the constraints containing $f(x, y) - f(x, z)$ by introducing a relaxation parameter. See \cite{Scholtes2001} for details about the relaxation method mentioned above for MPCC and see \cite{Li2023novel, Li2024solving} for details about the relaxation methods for WDP/MDP/eMDP.

Here, similarly as in  \cite{Li2023novel, Li2024solving}, we use
\begin{eqnarray*}\label{TWDP}
{\rm TWDP(t)}\qquad\min &&F(x, y)  \\
\mbox{s.t.} && (x,y)\in \Omega,~g(x, y) \le 0, ~h(x, y)=0,\\
&&f(x, y) - f(x, z) - u^T g(x, z) \leq t,~h(x, z) =0, \\
&&\nabla_z f(x, z) + \nabla_z g(x, z) u + \nabla_z h(x, z) v = 0, ~u\geq 0,
\end{eqnarray*}
\begin{eqnarray*}\label{TMDP}
{\rm TMDP(t)}\qquad\min &&F(x, y) \nonumber \\
\mbox{\rm s.t.} &&(x,y)\in \Omega,~g(x, y) \le 0, ~h(x, y)=0,\\
&&f(x, y) - f(x, z)\leq t,~ u^T g(x, z) \geq0,~h(x, z) =0,\nonumber\\
&&\nabla_z f(x, z) + \nabla_z g(x, z) u + \nabla_z h(x, z) v = 0, ~u\geq 0,\nonumber
\end{eqnarray*}
\begin{eqnarray*}\label{eTMDP}
{\rm eTMDP(t)}\qquad\min &&F(x, y) \nonumber \\
\mbox{\rm s.t.} &&(x,y)\in \Omega,~g(x, y) \le 0, ~h(x, y)=0,\\
&&f(x, y) - f(x, z)\leq t,~ u\circ g(x, z) \geq0, ~h(x, z) =0,\nonumber\\
&&\nabla_z f(x, z) + \nabla_z g(x, z) u + \nabla_z h(x, z) v = 0, ~u\geq 0\nonumber
\end{eqnarray*}
to approximate TWDP/TMDP/eTMDP respectively. Convergence results for these relaxation schemes are stated as follows.

\begin{thm}\label{1111}
Suppose that $\{t_k\}\downarrow0$ and $(x^k,y^k,z^k,u^k,v^k)$ is a globally optimal solution of ${\mathrm{TWDP}(t_k)}$/${\mathrm{TMDP}(t_k)}$/${\mathrm{eTMDP}(t_k)}$ for each $k$. Then, every accumulation point of $\{(x^k,y^k,z^k,u^k,v^k)\}$ is globally optimal to ${\rm TWDP/TMDP/eTMDP}$.
\end{thm}

\begin{thm}\label{thmRelaxation}
Suppose that  $\{t_k\}\downarrow0$ and $(x^k,y^k,z^k,u^k,v^k)$ be a KKT point of ${\mathrm{TWDP}(t_k)}$/${\mathrm{TMDP}(t_k)}$/${\mathrm{eTMDP}(t_k)}$ for each $k$.
Assume that $(\bar{x},\bar{y},\bar{z},\bar{u},\bar{v})$ is an accumulation point of $\{(x^k,y^k,z^k,u^k,v^k)\}$. If {\rm TWDP/TMDP/eTMDP} satisfies the constant positive linear dependence at $(\bar{x},\bar{y},\bar{z},\bar{u},\bar{v})$, then $(\bar{x},\bar{y},\bar{z},\bar{u},\bar{v})$ is a KKT point of ${\rm TWDP/TMDP/eTMDP}$.
\end{thm}

Theorem \ref{1111} is easy to verified. The proof of Theorem \ref{thmRelaxation} is provided in Appendix C.

The algorithm based on the above relaxation techniques with projection is described in Algorithm \ref{alg:One Projection Relaxation Algorithm}.

\begin{algorithm}
  \caption{\quad Relaxation Algorithm with Projection}
  \label{alg:One Projection Relaxation Algorithm}

\textbf{Step 0}  Choose an initial point $\tilde{x}^0\in X$, an initial relaxation parameter $t_0>0$, an update parameter $\sigma\in(0,1)$, and three accuracy tolerances $\{\epsilon_{\rm \scriptscriptstyle SQP},\epsilon_{\rm r},\epsilon_{\rm \scriptscriptstyle Inf}\}$. Set $k=0$.

\textbf{Step 1} 
Solve the lower-level program P$_{\tilde{x}^k}$ to get $(\tilde{y}^k,\tilde{u}^k,\tilde{v}^k)$. Set $\tilde{w}^k=(\tilde{x}^k,\tilde{y}^k,\tilde{y}^k,\tilde{u}^k,\tilde{v}^k)$.

\textbf{Step 2} Solve Q$(t_k)$ by \emph{fmincon} with $\tilde{w}^k$ as a starting point and $\epsilon_{\rm \scriptscriptstyle SQP}$ as an accuracy tolerance to obtain an iterative $w^k=(x^k,y^k,z^k,u^k,v^k)$.
\begin{itemize}
  \item If $w^k$ does not satisfy some termination criterion, go to Step 3.

  \item If $w^k$ satisfies the termination criterion, calculate {\tt Infeasibility} with $w^k$.
       If {\tt Infeasibility} $\leq \epsilon_{\rm \scriptscriptstyle Inf}$, stop. Otherwise, solve $x^*=\arg\min_{x \in X}\|x-x^k\|^2$ and solve the lower-level program P$_{x^*}$ to obtain a candidate solution $(x^*,y^*,y^*,u^*,v^*)$.
\end{itemize}

\textbf{Step 3}
Set $\tilde{x}^{k+1}={x}^k$ and $t_{k+1}=\max\{\sigma t_k, \epsilon_{\rm r}\}$. Return to Step 1 with $k=k+1$.
\end{algorithm}

Similarly to Algorithm \ref{alg:One Projection Direct Algorithm}, Algorithm \ref{alg:One Projection Relaxation Algorithm}
must output a feasible solution to the original BP as long as the lower-level programs are solvable and, in the case where some lower-level program has no solution, we need to change an initial point to try again.
Especially, in our experiments, the termination criterion in Step 2 was either $t_k\leq\epsilon_{\rm r}$ or
\begin{eqnarray*}
&& |f(x^k,y^k)- f(x^k, z^k) - {u^k}^T g(x^k, z^k)|\leq\epsilon_{\rm r}~~{\rm for~TWDP}(t_k), \\
&& |f(x^k,y^k)- f(x^k, z^k)|\leq\epsilon_{\rm r}~~{\rm for~TMDP}(t_k), \\
&& |f(x^k,y^k)- f(x^k, z^k)|\leq\epsilon_{\rm r}~~{\rm for~eTMDP}(t_k).
\end{eqnarray*}

To summarize, by adding the projection steps, both direct and relaxation algorithms can output feasible solutions to the original BP, which effectively addresses the shortcomings identified in the preceding studies  \cite{Li2023novel, Li2024solving}. Therefore, in our numerical comparison given in the next section, we will mainly focus on comparing the upper-level objective values and CPU times in evaluating different reformulations.

\section{Numerical comparison of seven single-level reformulations}
\label{sec-Numerical}
In this section, we compare seven single-level reformulations, including the old ones MPCC/WDP/MDP/eMDP and the new ones TWDP/TMDP/eTMDP by generating randomly various types of bilevel programs containing equality constraints. All experiments were implemented in MATLAB 9.13.0 and run on a computer with Windows 10, Intel(R) Xeon(R) E-2224G CPU $@$ 3.50GHz 3.50 GHz.

\subsection{Test problems}
We designed a procedure to generate randomly three types of bilevel programs as test problems, which include 150 bilevel programs with LP constraints, 150 bilevel programs with QP constraints, and 150 bilevel programs with QCQP (i.e., quadratically constrained quadratic program) constraints. All upper-level objectives and upper-level constraints are linear. In addition, all quadratic functions are convex.

\subsubsection{Group I: bilevel programs with LP constraints}
We generated 150 bilevel programs with LP constraints in the form
\begin{eqnarray}
\min &&c_1^T x+c_2^T y\nonumber\\
\mbox{s.t.}&&A_1x\leq b_1,\label{linear-linear BP}\\
&&y\in \arg\min\limits_{y}\{d_2^T y: A_2x+B_2y \leq b_2,A_3x+B_3y = b_3,~b_l\leq y\leq b_u\},\nonumber
\end{eqnarray}
where all coefficients $A_1\in \mathbb{R}^{l\times n}$, $A_2\in \mathbb{R}^{p\times n}$, $A_3\in \mathbb{R}^{q\times n}$, $B_2\in \mathbb{R}^{p\times m}$, $B_3\in \mathbb{R}^{q\times m}$, $b_1\in \mathbb{R}^l$, $b_2\in \mathbb{R}^p$, $b_3\in \mathbb{R}^q$, $c_1\in \mathbb{R}^n$, and $c_2,d_2\in \mathbb{R}^m$ were generated by \emph{sprand} in $[-1,1]$ with density to be 0.5, while the lower and upper bounds were generated by $b_l=-10*ones(m,1)$ and $b_u=10*ones(m,1)$ respectively.

The generated 150 examples contain three subgroups by adjusting only lower-level dimensions $\{m,p,q\}$, because these lower-level variables and constraints are the essential factors to cause difficulty in solving bilevel programs. The first 50 examples were generated with $\{n=20, l=25, m=100, p=110, q=20\}$, the second 50 examples were generated with $\{n=20, l=25, m=120, p=130, q=40\}$, and the last 50 examples were generated with $\{n=20, l=25, m=140, p=150, q=60\}$.

\subsubsection{Group II: bilevel programs with QP constraints}
We generated 150 bilevel programs with QP constraints in the form
\begin{eqnarray}
\min &&c_1^T x+c_2^T y\nonumber\\
\mbox{s.t.}&&A_1x\leq b_1,\label{linear-quadratic BP}\\
&&y\in \arg\min\limits_{y}\{\frac{1}{2}y^T H y+d_2^T y: A_2x+B_2y \leq b_2,A_3x+B_3y = b_3,~b_l\leq y\leq b_u\},\nonumber
\end{eqnarray}
where $H\in \mathbb{R}^{m\times m}$ is a symmetric positive semi-definite matrix and other coefficients are the same as in Case I. These 150 examples also contain three subgroups by adjusting lower-level dimensions only: The first 50 examples were generated with $\{n=20, l=25, m=30, p=20, q=10\}$, the second 50 examples were generated with $\{n=20, l=25, m=50, p=40, q=20\}$, and the last 50 examples were generated with $\{n=20, l=25, m=70, p=60, q=30\}$.

\subsubsection{Group III: bilevel programs with QCQP constraints}
We generated 150 bilevel programs with QCQP constraints in the form
\begin{eqnarray}
\min &&c_1^T x+c_2^T y\nonumber\\
\mbox{s.t.}&&A_1x\leq b_1,\label{linear-quadratic-quadratic BP}\\
&&y\in \arg\min\limits_{y}\{\frac{1}{2}y^T H y+d_2^T y:
A_2x+B_2y \leq b_2,A_3x+B_3y = b_3,~\frac{1}{2}y^T G y+d_4^T y\leq b_4,~b_l\leq y\leq b_u\},\nonumber
\end{eqnarray}
where $G,H\in \mathbb{R}^{m\times m}$ are symmetric positive semi-definite matrices, $d_4\in \mathbb{R}^m$, $b_4\in \mathbb{R}$, and other coefficients are the same as in Case I. The examples contain three subgroups by adjusting lower-level dimensions only: The first 50 examples were generated with $\{n=20, l=25, m=30, p=20, q=10\}$, the second 50 examples were generated with $\{n=20, l=25, m=50, p=40, q=20\}$, and the last 50 examples were generated with $\{n=20, l=25, m=70, p=60, q=30\}$.

\subsection{Numerical comparison}
In this subsection, we report our experiments on the algorithms given in Section 3 in solving the examples randomly generated in the above subsection. Our goal is to provide a comprehensive comparison of seven single-level reformulations for BP from multiple perspectives.

All numerical results for 450 test examples are shown in Tables 1-63 in the data repository materials, where `\#' represents the example serial number, `ObjVal' represents the upper-level objective value at the corresponding output, and `Time' represents the average CPU time (in seconds) by running three times to solve the same example. The accuracy tolerances in Algorithms \ref{alg:One Projection Direct Algorithm} and \ref{alg:One Projection Relaxation Algorithm} were set as {\tt Infeasibility} $<10^{-5}$. However, there were still a few test examples not to satisfy this accuracy tolerance and, 
for these examples, we weakened the accuracy tolerance as {\tt Infeasibility} $<10^{-4}$ so as to obtain satisfactory approximation solutions.

By adding the projection steps, except for a few test examples, Algorithms \ref{alg:One Projection Direct Algorithm} and \ref{alg:One Projection Relaxation Algorithm} could obtain approximation solutions feasible to the original bilevel programs.
Therefore, in our numerical comparison given below, when evaluating various algorithms, we mainly focus on the numbers of dominant cases and the CPU times. Here, `dominant case' means that the corresponding algorithm achieved either the best or one of the best approximation solutions by comparing their final objective values.

\subsubsection{Comparison of direct and relaxation algorithms with projection}
\label{subsubsec-Comparisondirectandrelaxation}
We first focus on comparing the direct and relaxation algorithms by applying Algorithms \ref{alg:One Projection Direct Algorithm} and \ref{alg:One Projection Relaxation Algorithm} to solve the numerical examples. The direct algorithm is obviously more time-saving and labor-saving than the relaxation algorithms. If its numerical performance is satisfactory, there is no reason to use the relaxation algorithms.

Comparison results of direct and relaxation algorithms with projection are summarized in Table \ref{tab:Comparison direct and relaxation}, where values in the
MPCC/WDP/MDP/eMDP/TWDP/TMDP/eTMDP columns represent the numbers of dominant cases of the corresponding algorithm. We do not compare the CPU times here because the direct algorithm is obviously more time-saving than the relaxation algorithms.

\begin{table}[htbp]
  \centering
  \caption{Comparison of direct and relaxation algorithms with projection}
    \tabcolsep=1mm
    \begin{threeparttable}
    \begin{tabular}{cccccccccc}
    \hline\specialrule{0.04em}{2pt}{4pt}
 Group & Dimension  & Method & MPCC & WDP & MDP & eMDP & TWDP & TMDP & eTMDP\\
    \specialrule{0.04em}{3pt}{2pt}
    \multirow{6}{*}{ Group I}
   & \multirow{2}[0]{*}{$m=100, p=110, q=20$}
      & Direct     & 18  & 0   & 1   & 0   & 1   & 0  & 0 \\
   &  & Relaxation & 32  & 50  & 49  & 50  & 49  & 50 & 50\\

   & \multirow{2}[0]{*}{$m=120, p=130, q=40$}
      & Direct     & 21  & 0   & 0   & 0   & 0   & 1  & 0 \\
   &  & Relaxation & 29  & 50  & 50  & 50  & 50  & 49 & 50\\

   & \multirow{2}[0]{*}{$m=140, p=150, q=60$~\!\tnote{\tt1}}
      & Direct     & 17  & 0   & 0   & 0   & 0   & 0  & 0 \\
   &  & Relaxation & 31  & 48  & 48  & 48  & 48  & 48 & 48\\

   \cmidrule(r){3-10}
    \multirow{6}{*}{ Group II}
   & \multirow{2}[0]{*}{$m=30, p=20, q=10$~\!\tnote{\tt2}}
      & Direct     & 9   & 2   & 1   & 1   & 1   & 1  & 0\\
   &  & Relaxation & 41  & 48  & 50  & 50  & 49  & 49 & 50\\

   & \multirow{2}[0]{*}{$m=50, p=40, q=20$}
      & Direct     & 12  & 2   & 2   & 1   & 1   & 0  & 1 \\
   &  & Relaxation & 38  & 48  & 48  & 49  & 49  & 50 & 49\\

   & \multirow{2}[0]{*}{$m=70, p=60, q=30$}
      & Direct     & 25  & 0   & 0   & 4   & 0   & 1  & 1 \\
   &  & Relaxation & 25  & 50  & 50  & 47  & 50  & 49 & 49\\

   \cmidrule(r){3-10}
    \multirow{6}{*}{ Group III}
   & \multirow{2}[0]{*}{$m=30, p=20, q=10$}
      & Direct     & 23  & 5   & 3   & 30   & 6   & 3 & 22\\
   &  & Relaxation & 27  & 46  & 47  & 20   & 45  & 47& 28\\

   & \multirow{2}[0]{*}{$m=50, p=40, q=20$}
      & Direct     &  19 & 7   & 5   & 30   & 9   & 17& 33\\
   &  & Relaxation &  31 & 43  & 47  & 20   & 41  & 33& 17\\

   & \multirow{2}[0]{*}{$m=70, p=60, q=30$}
      & Direct     & 29  & 12  & 9   & 24   & 12  & 13& 22\\
   &  & Relaxation &  21 & 39  & 41  & 26   & 39  & 38& 28\\
   \hline\hline
    \multirow{2}{*}{Totals}
      & \multirow{2}[0]{*}{}
      & Direct     & 173  & 28  & 21   & 90   & 30  & 36& 79\\
   &  & Relaxation &  275 & 422  & 430  & 360   & 420  & 413& 369\\
    \hline\hline
    \end{tabular}%
  \label{tab:Comparison direct and relaxation}
\begin{tablenotes}[flushleft]
\footnotesize
\item[\tt1] The total dominant cases less than the total number 50 means that some examples are infeasible.  
\item[\tt2] The total dominant cases larger than the total number 50 means that, for some examples, both direct and relaxation algorithms generated the best approximation solutions. 
\end{tablenotes}
\end{threeparttable}

\end{table}

It can be seen from Table \ref{tab:Comparison direct and relaxation} that, at least in our experiments, the relaxation algorithms were significantly better than the direct algorithm from an overall perspective, especially in the cases of duality-based reformulations. Among these duality-based reformulations, the extended two ones eMDP and eTMDP performed weaker than the others. This phenomenon coincides with our observations given in the next two subsection and contradicts our initial conjecture observed from the numerical experiments in \cite{Li2024solving} that the tighter the feasible region of the reformulation, the better the corresponding numerical performance. There are two possible reasons: One is that introducing the projection steps may cause some changes in numerical performances of different reformulations. The other is that limited examples given in \cite{Li2024solving} may lead to some deviations in evaluating the algorithms.

\subsubsection{Overall comparison of seven single-level reformulations}

In this subsection, we compare the overall numerical performances of seven single-level reformulations and, in the next subsection, we give a pairwise comparison of them. Overall comparison results of seven single-level reformulations are summarized in Tables \ref{tab:Comparison existing reformulations only method Case I}-\ref{tab:Comparison existing reformulations only method totals}.

\begin{table}[htbp]
  \centering
  \tabcolsep=1.48mm
  \caption{Overall comparison of seven reformulations for direct or relaxation algorithms in Group I}
    \begin{tabular}{ccccccc}
    \hline\hline
    \multirow{2}{*}{Dimension}&\multirow{2}{*}{Reformulation}& {Dominant cases} & Average time (s)& {Dominant cases} & Average time (s) \\
    &&(direct) & (direct) &(relaxation) &(relaxation)& \\
    \hline
    \multirow{7}[2]{*}{$\begin{matrix}m=100\\ p=110\\q=20\end{matrix}$}
    & MPCC  & 2   & 10.15  &  2    & 427.69  \\
    & WDP   & 10  & 13.25  &  14   & 205.49  \\
    & MDP   & 8   & 16.63  & 16   & 171.81  \\
    & eMDP  & 3   & 78.78  &  10   & 530.86  \\
    & TWDP  & 15  & 16.86  &  12   & 212.90  \\
    & TMDP  & 11  & 15.62  &  18   & 206.32  \\
    & eTMDP & 4   & 96.19  &  5    & 469.42  \\

    \cmidrule(r){2-4} \cmidrule(r){5-7}
    \multirow{7}[2]{*}{$\begin{matrix}m=120\\ p=130\\ q=40\end{matrix}$}
    & MPCC  & 5   & 10.49  &  1    & 502.13  \\
    & WDP   & 9   & 17.79  &  12   & 260.57  \\
    & MDP   & 11  & 26.34  &  13   & 325.02  \\
    & eMDP  & 3   & 201.96 &  12   & 1356.16 \\
    & TWDP  & 18  & 22.49  &  21   & 304.25  \\
    & TMDP  & 12  & 30.74  &  10   & 397.22  \\
    & eTMDP & 6   & 177.14 &  9    & 1265.11 \\

    \cmidrule(r){2-4} \cmidrule(r){5-7}
    \multirow{7}[2]{*}{$\begin{matrix}m=140\\p=150\\ q=60\end{matrix}$}
    & MPCC  & 7   & 19.14  &  3    & 586.87  \\
    & WDP   & 14  & 47.78  &  17   & 518.32  \\
    & MDP   & 6   & 35.78  &  16   & 430.77  \\
    & eMDP  & 8   & 1174.12&  9    & 2576.61 \\
    & TWDP  & 10  & 57.57  &  9    & 400.55  \\
    & TMDP  & 11  & 55.12  &  13   & 428.28  \\
    & eTMDP & 8   & 869.14 &  7    & 3073.16 \\
    \hline\hline
    \end{tabular}
  \label{tab:Comparison existing reformulations only method Case I}
\end{table}

\begin{table}[htbp]
  \centering
  \tabcolsep=1.7mm
  \caption{Overall comparison of seven reformulations for direct or relaxation algorithms in Group II}
  \begin{tabular}{ccccccc}
   \hline\hline
   \multirow{2}{*}{Dimension}&\multirow{2}{*}{Reformulation}& {Dominant cases} & Average time (s)& {Dominant cases} & Average time (s) \\
    &&(direct) & (direct) &(relaxation) &(relaxation)& \\
   \hline

    \multirow{7}[0]{*}{$\begin{matrix}m=30\\ p=20\\ q=10\end{matrix}$}
    & MPCC  & 10    & 0.88  &  11   & 40.61  \\
    & WDP   & 10    & 1.87  &  20   & 24.90  \\
    & MDP   & 9     & 1.92  &  29   & 16.34  \\
    & eMDP  & 4     & 2.36  &  12   & 23.87  \\
    & TWDP  & 6     & 1.83  &  27   & 24.80  \\
    & TMDP  & 9     & 1.92  &  23   & 15.06  \\
    & eTMDP & 4     & 2.40  &  10   & 23.37  \\

    \cmidrule(r){2-4} \cmidrule(r){5-7}
    \multirow{7}[0]{*}{$\begin{matrix}m=50\\ p=40\\ q=20\end{matrix}$}
    & MPCC  & 7     & 1.88  &  7    & 119.62 \\
    & WDP   & 10    & 5.96  &  18   & 98.33  \\
    & MDP   & 10    & 6.39  &  13   & 97.31  \\
    & eMDP  & 5     & 9.56  &  3    & 198.06 \\
    & TWDP  & 5     & 6.23  &  9    & 93.64  \\
    & TMDP  & 13    & 11.59 &  11   & 84.20  \\
    & eTMDP & 2     & 12.71 &  8    & 211.74 \\

    \cmidrule(r){2-4} \cmidrule(r){5-7}
    \multirow{7}[0]{*}{$\begin{matrix}m=70\\ p=60\\ q=30\end{matrix}$}
    & MPCC  & 7     & 3.29  &  2    & 235.52 \\
    & WDP   & 11    & 15.04 &  15   & 293.94 \\
    & MDP   & 14    & 24.10 &  16   & 330.52 \\
    & eMDP  & 7     & 49.54 &  9    & 734.36 \\
    & TWDP  & 6     & 15.23 &  17   & 338.77 \\
    & TMDP  & 15    & 18.06 &  10   & 353.84 \\
    & eTMDP & 3     & 42.12 &  11   & 902.82 \\
    \hline\hline
    \end{tabular}
  \label{tab:Comparison existing reformulations only method Case II}
\end{table}%

\begin{table}[htbp]
  \centering
  \tabcolsep=1.7mm
  \caption{Overall comparison of seven reformulations for direct or relaxation algorithms in Group III}
  \begin{tabular}{ccccccc}
   \hline\hline
  \multirow{2}{*}{Dimension}&\multirow{2}{*}{Reformulation}& {Dominant cases} & Average time (s)& {Dominant cases} & Average time (s) \\
    &&(direct) & (direct) &(relaxation) &(relaxation)& \\
   \hline

   \multirow{7}[0]{*}{$\begin{matrix}m=30\\ p=20\\ q=10\end{matrix}$}
   & MPCC  & 7     & 1.50  &  11   & 94.79  \\
   & WDP   & 5     & 2.78  &  20   & 142.69  \\
   & MDP   & 3     & 2.63  &  14   & 128.53  \\
   & eMDP  & 19    & 3.14  &  1    & 81.55  \\
   & TWDP  & 3     & 2.40  &  7    & 144.74  \\
   & TMDP  & 7     & 2.56  &  6    & 136.54  \\
   & eTMDP & 20    & 2.67  &  1    & 70.45  \\

   \cmidrule(r){2-4} \cmidrule(r){5-7}
   \multirow{7}[0]{*}{$\begin{matrix}m=50\\ p=40\\ q=20\end{matrix}$}
   & MPCC  & 7     & 3.41  &  15   & 175.00  \\
   & WDP   & 1     & 6.43  &  18   & 403.96  \\
   & MDP   & 0     & 9.72  &  12   & 282.87  \\
   & eMDP  & 19    & 16.36 &  0    & 526.95  \\
   & TWDP  & 6     & 6.87  &  7    & 365.16  \\
   & TMDP  & 5     & 8.41  &  15   & 201.37  \\
   & eTMDP & 21    & 9.27  &  2    & 357.25  \\

   \cmidrule(r){2-4} \cmidrule(r){5-7}
   \multirow{7}[0]{*}{$\begin{matrix}m=70\\ p=60\\ q=30\end{matrix}$}
   & MPCC  & 11    & 8.91  &  10   & 343.72  \\
   & WDP   & 6     & 19.69 &  19   & 889.89  \\
   & MDP   & 6     & 17.54 &  10   & 793.53  \\
   & eMDP  & 11    & 50.18 &  0    & 1727.30 \\
   & TWDP  & 3     & 14.97 &  3    & 828.71  \\
   & TMDP  & 5     & 22.43 &  9    & 774.16  \\
   & eTMDP & 12    & 36.43 &  4    & 2377.36 \\
   \hline\hline
   \end{tabular}
  \label{tab:Comparison existing reformulations only method Case III}
\end{table}

\begin{table}[htbp]
  \centering
  \tabcolsep=1.7mm
\begin{threeparttable}
  \caption{Overall comparison of total dominant cases of seven reformulations}
  \begin{tabular}{ccccccccc}
   \hline\hline
&&\multirow{2}{*}{Reformulation}&& {Total dominant cases} && {Total dominant cases}  &&\\
& &   && (direct) &&(relaxation) &&\\
   \hline
 && MPCC  && 63     && 62    &&\\
 &&    WDP   && 76     && 153   && \\
  &&   MDP   && 67     && 139  &&   \\
 &&    eMDP  && 79    && 56    && \\
 &&    TWDP  && 72    && 112  &&  \\
 &&    TMDP  && 88     && 115  &&  \\
 &&    eTMDP && 80    && 57  &&  \\
   \hline\hline
   \end{tabular}
  \label{tab:Comparison existing reformulations only method totals}
\begin{tablenotes}[flushleft]
\footnotesize
\item[\tt1] The total dominant cases larger than the total number 450 means that, for some examples, more than one reformulation generated the best approximation solutions.
\end{tablenotes}
\end{threeparttable}
\end{table}

We first analyze the numerical performance of direct algorithm. In terms of dominant data, at least in our experiments, although there was no significant difference in their total dominant cases from Table \ref{tab:Comparison existing reformulations only method totals}, the duality-based reformulations performed slightly better than the MPCC reformulation. However, the situation in terms of CPU times was quite different. It can be seen from Tables \ref{tab:Comparison existing reformulations only method Case I}-\ref{tab:Comparison existing reformulations only method Case III} that seven reformulations had consistent results in solving various subgroups, that is, MPCC was always the best one, while the extended reformulations eMDP and eTMDP were the worst ones.

We next analyze the numerical performance of relaxation algorithms.
It can be seen from Table \ref{tab:Comparison existing reformulations only method totals} that, in terms of dominant data, the WDP/MDP/TWDP/TMDP reformulations were more than twice as good as the MPCC/eMDP/eTMDP reformulations, but the differences among WDP/MDP/TWDP/TMDP or MPCC/eMDP/eTMDP are not significant.
On the other hand, Tables \ref{tab:Comparison existing reformulations only method Case I}-\ref{tab:Comparison existing reformulations only method Case III} reveal that, in terms of CPU times, the WDP/MDP/TWDP/TMDP reformulations were also better than the MPCC reformulation, while the latter was more time-saving than the eMDP/eTMDP reformulations and the gaps among them were getting bigger and bigger with the increase of data scales.

Note that the extended two reformulations eMDP and eTMDP are always the worst ones among the duality-based reformulations so far.

\subsubsection{Pairwise comparison of seven single-level reformulations}
In order to compare seven reformulations in more depth, we provide a pairwise comparison of them in this subsection.
Pairwise comparison results of seven single-level reformulations for direct and relaxation algorithms
are shown in Tables \ref{tab:Pairwise comparison direct}-\ref{tab:Pairwise comparison all cases-relaxation}, where
`$\frac{\mathrm{A}}{\mathrm{B}}$' denotes the quotient of two dominant data of the A and B reformulations. In particular,
Tables \ref{tab:Pairwise comparison direct}-\ref{tab:Pairwise comparison relaxation} provide the detailed data in different groups, while Tables \ref{tab:Pairwise comparison all cases-direct}-\ref{tab:Pairwise comparison all cases-relaxation} summarize the total data for 450 numerical examples, where a number more than one means that A outperformed B in our experiments and conversely,  a number less than one means that B outperformed A in our experiments.

\begin{table}[htbp]
  \centering
  \tabcolsep=2.4mm
  \caption{Pairwise comparison of dominant cases ($\frac{\mathrm{A}}{\mathrm{B}}$) for direct algorithm in Groups I/II/III}
    \begin{tabular}{cccccccc}
     \hline\hline
    \diagbox{A}{B}
   & MPCC & WDP & MDP& eMDP & TWDP& TMDP \\
     \specialrule{0.04em}{0pt}{4pt}

WDP   & $\frac{105}{67}$/$\frac{90}{64}$/$\frac{71}{79}$
           & -   &  -  &   - &  -  &  -    \\[2mm]

MDP   & $\frac{101}{66}$/$\frac{94}{62}$/$\frac{73}{78}$

      & $\frac{90}{91}$/$\frac{77}{78}$/$\frac{77}{78}$
                 & -   &  -  &  -  &  -    \\[2mm]

eMDP  & $\frac{55}{94}$/$\frac{82}{70}$/$\frac{77}{73}$
      & $\frac{42}{107}$/$\frac{53}{99}$/$\frac{78}{72}$
      & $\frac{44}{106}$/$\frac{56}{96}$/$\frac{81}{69}$
                       & -   &  -  &  -   \\[2mm]

TWDP  & $\frac{118}{57}$/$\frac{92}{64}$/$\frac{78}{73}$
      & $\frac{99}{81}$/$\frac{81}{76}$/$\frac{80}{77}$
      & $\frac{99}{81}$/$\frac{61}{97}$/$\frac{80}{76}$
      & $\frac{109}{40}$/$\frac{95}{57}$/$\frac{74}{76}$
                             & -   &  -   \\[2mm]

TMDP  & $\frac{110}{63}$/$\frac{98}{58}$/$\frac{80}{71}$
      & $\frac{95}{89}$/$\frac{89}{72}$/$\frac{85}{69}$
      & $\frac{99}{86}$/$\frac{78}{76}$/$\frac{81}{74}$
      & $\frac{107}{42}$/$\frac{100}{52}$/$\frac{70}{81}$
      & $\frac{85}{92}$/$\frac{82}{76}$/$\frac{84}{73}$
                                   & -    \\[2mm]

eTMDP & $\frac{61}{88}$/$\frac{73}{80}$/$\frac{82}{68}$
      & $\frac{55}{94}$/$\frac{44}{107}$/$\frac{80}{70}$
      & $\frac{52}{98}$/$\frac{42}{110}$/$\frac{84}{66}$
      & $\frac{89}{72}$/$\frac{66}{85}$/$\frac{100}{90}$
      & $\frac{41}{109}$/$\frac{53}{99}$/$\frac{78}{72}$
      & $\frac{47}{103}$/$\frac{42}{109}$/$\frac{85}{65}$
                                          \\[1mm]
    \hline\hline
    \end{tabular}
  \label{tab:Pairwise comparison direct}
\end{table}

\begin{table}[htbp]
  \centering
  \tabcolsep=1.6mm
  \caption{Pairwise comparison of dominant cases ($\frac{\mathrm{A}}{\mathrm{B}}$) for relaxation algorithms in Groups I/II/III}
    \begin{tabular}{cccccccc}
     \hline\hline
    \diagbox{A}{B}
   & MPCC & WDP & MDP& eMDP & TWDP& TMDP  \\
     \specialrule{0.04em}{0pt}{4pt}

WDP   & $\frac{133}{15}$/$\frac{128}{26}$/$\frac{112}{38}$
           & -   &  -  &   - &  -  &  -    \\[2mm]

MDP   & $\frac{139}{9}$/$\frac{126}{29}$/$\frac{111}{39}$
      & $\frac{100}{86}$/$\frac{104}{97}$/$\frac{71}{95}$
                 & -   &  -  &  -  &  -    \\[2mm]

eMDP  & $\frac{132}{16}$/$\frac{114}{37}$/$\frac{58}{92}$
      & $\frac{60}{98}$/$\frac{57}{113}$/$\frac{7}{143}$
      & $\frac{50}{108}$/$\frac{58}{112}$/$\frac{4}{146}$
                       & -   &  -  &  -   \\[2mm]

TWDP  & $\frac{136}{12}$/$\frac{127}{28}$/$\frac{101}{49}$
      & $\frac{91}{94}$/$\frac{106}{98}$/$\frac{53}{109}$
      & $\frac{92}{90}$/$\frac{105}{100}$/$\frac{54}{106}$
      & $\frac{103}{51}$/$\frac{106}{61}$/$\frac{142}{8}$
                             & -   &  -   \\[2mm]

TMDP  & $\frac{135}{13}$/$\frac{129}{25}$/$\frac{124}{76}$
      & $\frac{98}{92}$/$\frac{101}{105}$/$\frac{50}{106}$
      & $\frac{93}{90}$/$\frac{99}{102}$/$\frac{56}{98}$
      & $\frac{99}{55}$/$\frac{112}{59}$/$\frac{131}{19}$
      & $\frac{100}{90}$/$\frac{103}{105}$/$\frac{69}{87}$
                                   & -     \\[2mm]

eTMDP & $\frac{127}{21}$/$\frac{115}{36}$/$\frac{68}{82}$
      & $\frac{60}{100}$/$\frac{62}{108}$/$\frac{11}{139}$
      & $\frac{49}{110}$/$\frac{63}{110}$/$\frac{12}{138}$
      & $\frac{79}{86}$/$\frac{86}{87}$/$\frac{101}{72}$
      & $\frac{54}{102}$/$\frac{67}{108}$/$\frac{19}{135}$
      & $\frac{60}{96}$/$\frac{65}{104}$/$\frac{24}{127}$
                                          \\[1mm]
    \hline\hline
    \end{tabular}
  \label{tab:Pairwise comparison relaxation}
\end{table}

%
%
%
%
%
%

\begin{table}
  \centering
  \tabcolsep=3.2mm
  \caption{Pairwise comparison of ratios  ($\frac{\mathrm{A}}{\mathrm{B}}$) for direct algorithm with 450 examples}
    \begin{tabular}{cccccccc}
     \hline\hline
    \diagbox{A}{B}
   & MPCC & WDP & MDP& eMDP & TWDP& TMDP  \\
     \specialrule{0.04em}{0pt}{4pt}

WDP   & 1.27 & - &  -   & -   &   -   & -    \\[2mm]

MDP   & 1.30 & 0.99 & - & -  &  - &  -   \\[2mm]

eMDP  & 0.90 & 0.62 & 0.67 &- &-&- \\[2mm]

TWDP  & 1.48 & 1.11 & 0.95 & 1.61
      &-&- \\[2mm]

TMDP  & 1.50 & 1.17 & 1.09 & 1.58
      & 1.04 & -   \\[2mm]

eTMDP &0.92 &0.66 &0.65 &1.03 &0.61
      &0.63  \\[1mm]
    \hline\hline
    \end{tabular}
  \label{tab:Pairwise comparison all cases-direct}
\end{table}

\begin{table}
  \centering
  \tabcolsep=3.2mm
  \caption{Pairwise comparison of ratios  ($\frac{\mathrm{A}}{\mathrm{B}}$) for relaxation algorithm with 450 examples}
    \begin{tabular}{cccccccc}
     \hline\hline
    \diagbox{A}{B}
   & MPCC & WDP & MDP& eMDP & TWDP& TMDP  \\
     \specialrule{0.04em}{0pt}{4pt}

WDP   & 4.72 & - &  -   & -   &   -   & -    \\[2mm]

MDP   & 4.77 & 0.99 & - & -  &  - &  -   \\[2mm]

eMDP  & 2.10 & 0.35 & 0.31 &- &-&- \\[2mm]

TWDP  & 4.09 & 0.83 & 0.85 & 2.93
      &-&- \\[2mm]

TMDP  & 3.40 & 0.82 & 0.86 & 2.57
      & 0.96 & -   \\[2mm]

eTMDP &2.23 &0.38 &0.35 &1.09 &0.41
      &0.46  \\[1mm]
    \hline\hline
    \end{tabular}
  \label{tab:Pairwise comparison all cases-relaxation}
\end{table}

Table \ref{tab:Pairwise comparison all cases-direct} reports a pairwise comparison for direct algorithm with 450 examples. It can be seen that, in our experiments, the WDP/MDP/TWDP/TMDP reformulations outperformed the MPCC reformulation, while the latter outperformed the eMDP/eTMDP reformulations.
The relationship chain among seven reformulations was as follows:
$$
{\rm TMDP} \ \succ \ {\rm TWDP} \ \succ \ {\rm WDP} \ \succ \ {\rm MDP} \ \succ \ {\rm MPCC} \ \succ \ {\rm eTMDP} \ \succ \ {\rm eMDP},
$$
where ${\rm A}\succ {\rm B}$ means that A outperformed B in our experiments.

Table \ref{tab:Pairwise comparison all cases-relaxation} gives a pairwise comparison for relaxation algorithms with 450 examples. It can be seen that, in our experiments, all duality-based reformulations outperformed the MPCC reformulation. In particular,
the WDP/MDP/TWDP/TMDP reformulations were 3-5 times better than the MPCC reformulation, while the eMDP/eTMDP reformulations were 2 times better than the MPCC reformulation.
The relationship chain among seven reformulations was as follows:
$$
{\rm WDP} \ \succ \ {\rm MDP} \ \succ \ {\rm TWDP} \ \succ \ {\rm TMDP} \ \succ \ {\rm eTMDP} \ \succ \ {\rm eMDP} \ \succ \ {\rm MPCC}.
$$

To sum up, whether overall comparison or pairwise comparison, and whether direct algorithm or relaxation algorithms, at least in our experiments, the WDP/MDP/TWDP/TMDP reformulations were always better than the MPCC reformulation, while the extended two reformulations eMDP and eTMDP were always the worst ones among the duality-based reformulations, which shows that our conjecture observed from the numerical experiments in \cite{Li2024solving} is incorrect, that is, the tighter the feasible region of the reformulation, the better the corresponding numerical performance. One reason may be that introducing the projection steps in  Algorithms \ref{alg:One Projection Direct Algorithm} and \ref{alg:One Projection Relaxation Algorithm} causes some changes in numerical performances of different reformulations.

\section{Conclusions}
\label{sec-conclusion}
Inspired by a conjecture from \cite{Li2024solving} that the tighter feasible region of a reformulation, the better its numerical performance, we proposed three new single-level reformulations TWDP/TMDP/eTMDP with tighter feasible regions by moving the lower-level equality constraints outside from WDP/MDP/eMDP respectively. These new ones inherit one main advantage of WDP/MDP/eMDP, that is, unlike MPCC, they may satisfy the MFCQ at their feasible points. 

Our main goal is to compare seven single-level reformulations, including the old ones MPCC/WDP/MDP/eMDP and the new ones TWDP/TMDP/eTMDP by numerical experiments. Different from the algorithms in \cite{Li2023novel, Li2024solving}, we add some projection steps into the algorithms as an improvement so as to ensure to output feasible points of BP always in theory at least.
Through numerical experiments on 450 test examples generated randomly, we show that, whether overall comparison or pairwise comparison, at least in our tests, the WDP/MDP/TWDP/TMDP reformulations were always better than the MPCC reformulation, while the extended two reformulations eMDP and eTMDP were always the worst ones among the duality-based reformulations, which shows that our conjecture observed from the numerical experiments in \cite{Li2024solving} is incorrect. In particular, for the relaxation algorithms, the WDP/MDP/TWDP/TMDP reformulations were 3-5 times better than the MPCC reformulation, while the eMDP/eTMDP reformulations were 2 times better than the MPCC reformulation.
In the next step, we will try to apply the duality-based approach to more general bilevel programs like pessimistic bilivel programs and mixed-integer bilevel programs, or other optimization problems closely to bilevel programs such as minimax problems and two-stage stochastic programming problems.


\section*{Acknowledgment}
This work was supported in part by 
the China Postdoctoral Science Foundation (No. 2024M761920).

\begin{appendices}

\setcounter{equation}{0}
\renewcommand\theequation{A.\arabic{equation}}
\setcounter{thm}{0}
\renewcommand\thethm{A.\arabic{thm}}
\setcounter{defi}{0}
\renewcommand\thedefi{A.\arabic{defi}}
\setcounter{rem}{0}
\renewcommand\therem{A.\arabic{rem}}

\section*{Appendix A: Proofs of Theorems 2.1--2.3}

\noindent {\bf Proof of Theorem 2.1}
We first show the weak duality. Since $\nabla_z L(x,z, u, v)=0$ for any $y\in Y(x)$ and $(z,u,v)\in \mathcal{F}_1(x)$,
it follows from the pseudoconvexity of $L(x,\cdot, u, v)$ that $L(x,y, u, v)\geq L(x,z, u, v)$.
Then, for $u\in \mathbb{R}^p_+$, $y\in Y(x)$,  and $(z,u,v)\in\mathcal{F}_1(x)$, we have
\begin{eqnarray*}
f(x,y)\geq f(x,y) + u^T g(x,y)= L(x,y, u, v)\geq L(x,z, u, v)=f(x,z) + u^T g(x,z) .
\end{eqnarray*}
By the arbitrariness of $y\in Y(x)$ and $(z,u,v)\in\mathcal{F}_1(x)$, we have
\begin{eqnarray*}\label{mWD-dualitythm-3}
\min\limits_{y\in Y(x)} f(x,y)\geq \max\limits_{(z, u, v)\in \mathcal{F}_1(x)} f(x,z) + u^T g(x,z).
\end{eqnarray*}

Next, we show the strong duality.
Since $\mathrm{P}_x$ satisfies the GCQ at $y_x \in S(x)$,
there is $(u_x,v_x)$ such that
\begin{eqnarray*}
\nabla_y L(x,y_x,u_x,v_x)=0,~u_x^T g(x,y_x)=0, ~u_x\geq 0,~h(x,y_x)=0,
\end{eqnarray*}
from which we have $(y_x,u_x,v_x)\in \mathcal{F}_1(x)$. Then, we have
\begin{eqnarray*}
\min\limits_{y\in Y(x)} f(x,y)\geq \max\limits_{(z, u, v)\in \mathcal{F}_1(x)} f(x,z) + u^T g(x,z)\geq f(x,y_x) + u_x^T g(x,y_x)=f(x,y_x).
\end{eqnarray*}
Since $y_x\in S(x)$, the above inequalities are all equalities. Therefore, we have
\begin{eqnarray*}
\min\limits_{y\in Y(x)} f(x,y)=f(x,y_x)= f(x,y_x) + u_x^T g(x,y_x) = \max\limits_{(z, u, v)\in \mathcal{F}_1(x)} f(x,z) + u^T g(x,z).
\end{eqnarray*}
This completes the proof.

%

\noindent {\bf Proof of Theorem 2.2}
We first show the weak duality. For any $y\in Y(x)$ and $(z,u,v)\in \mathcal{F}_2(x)$, we have
\begin{eqnarray}\label{mMD-dualitythm-1}
u^T g(x,y) + v^T h(x,y)-(u^T g(x,z) + v^T h(x,z))\leq0.
\end{eqnarray}
Since $u^Tg(x,\cdot)+v^Th(x,\cdot)$ is quasiconvex, \eqref{mMD-dualitythm-1} implies
\begin{eqnarray}\label{mMD-dualitythm-2}
(y-z)^T(\nabla g(x,z)u + \nabla h(x,z)v)\leq0.
\end{eqnarray}
Since $\nabla f(x,z)+\nabla g(x,z)u + \nabla h(x,z)v=0$, it follows that
\begin{eqnarray}\label{mMD-dualitythm-3}
(y-z)^T(\nabla f(x,z)+\nabla g(x,z)u + \nabla h(x,z)v)=0.
\end{eqnarray}
\eqref{mMD-dualitythm-2}-\eqref{mMD-dualitythm-3} implies
$(y-z)^T\nabla f(x,z)=(y-z)^T(\nabla g(x,z)u + \nabla h(x,z)v)\geq0.$
Since $f(x,\cdot)$ is pseudoconvex, we have $f(x,y)\geq f(x,z)$. By the arbitrariness of $y\in Y(x)$ and $(z,u,v)\in\mathcal{F}_2(x)$, we have
\begin{eqnarray*}
\min\limits_{y\in Y(x)} f(x,y)\geq \max\limits_{(z, u, v)\in \mathcal{F}_2(x)} f(x,z).
\end{eqnarray*}

Next, we show the strong duality.
Since $\mathrm{P}_x$ satisfies the GCQ at $y_x\in S(x)$,
there is $(u_x,v_x)$ such that
$\nabla_y L(x,y_x,u_x,v_x)=0,~u_x^T g(x,y_x)=0, ~u_x\geq 0,~h(x,y_x)=0,$
from which we have $(y_x,u_x,v_x)\in \mathcal{F}_2(x)$. Then, we have
\begin{eqnarray*}
\min\limits_{y\in Y(x)} f(x,y)\geq \max\limits_{(z, u, v)\in \mathcal{F}_2(x)} f(x,z)\geq f(x,y_x) .
\end{eqnarray*}
Since $y_x\in S(x)$, the above inequalities are all equalities. Therefore, we have
\begin{eqnarray*}
\min\limits_{y\in Y(x)} f(x,y)=f(x,y_x)= \max\limits_{(z, u, v)\in \mathcal{F}_2(x)} f(x,z).
\end{eqnarray*}
This completes the proof.

%
%

\noindent {\bf Proof of Theorem 2.3} Since the proof of this theorem is similar to Theorem 2.2 and so we omit it here.

\setcounter{equation}{0}
\renewcommand\theequation{B.\arabic{equation}}
\setcounter{thm}{0}
\renewcommand\thethm{B.\arabic{thm}}
\setcounter{defi}{0}
\renewcommand\thedefi{B.\arabic{defi}}
\setcounter{table}{0}
\renewcommand\thetable{B.\arabic{table}}
\setcounter{ex}{0}
\renewcommand\theex{B.\arabic{ex}}
\section*{Appendix B: Proofs of Theorems 3.2 and 3.3}

 For simplicity, we ignore the upper constraint $(x,y)\in \Omega$ throughout this appendix.

\noindent {\bf Proof of Theorem 3.2}
(i) Note that the case of MDP has been shown in Theorem 9 of Li et al. (2024), 
while other cases of eMDP/TWDP/TMDP/eTMDP can be shown in a similar way. For completeness, we provide a proof for the case of eMDP as an example.
In fact, if $(\bar{x},\bar{y},\bar{z},\bar{u},\bar{v})$ is a KKT point of WDP, there is $(\zeta^g, \zeta^h, \zeta^u, \widehat{\alpha}, \widehat{\beta}) \in \mathbb{R}^{2p+q+1+m}$ such that
\begin{eqnarray}
\label{WDPKKT-1}&&\nabla_x  F(\bar{x},\bar{y})  +\widehat{\alpha}(\nabla_x  f(\bar{x},\bar{y})-\nabla_x L(\bar{x},\bar{z}, \bar{u},\bar{v}))+ \nabla_{zx}^2 L(\bar{x},\bar{z}, \bar{u},\bar{v}) \widehat{\beta} +\nabla_x g(\bar{x},\bar{y}) \zeta^g+\nabla_x h(\bar{x},\bar{y}) \zeta^h=0,\qquad\\
\label{WDPKKT-2}&&\nabla_y F(\bar{x},\bar{y}) + \widehat{\alpha}\nabla_y f(\bar{x},\bar{y}) + \nabla_y g(\bar{x},\bar{y})\zeta^g+\nabla_y h(\bar{x},\bar{y})\zeta^h=0,\\
\label{WDPKKT-3}&&\nabla_{zz}^2 L(\bar{x},\bar{z}, \bar{u},\bar{v}) \widehat{\beta}=0, \\
\label{WDPKKT-4}&&-\hat{\alpha} g(\bar{x},\bar{z}) + \nabla_z g(\bar{x},\bar{z})^T\widehat{\beta} - \zeta^u=0, \\
\label{WDPKKT-5}&&-\widehat{\alpha} h(\bar{x},\bar{z})+\nabla_z h(\bar{x},\bar{z})^T\widehat{\beta} =0,\\
\label{WDPKKT-6}&&0\le \widehat{\alpha} ~\bot ~f(\bar{x},\bar{y})-L(\bar{x},\bar{z},\bar{u},\bar{v}) \leq0,\\
\label{WDPKKT-7}&&0\le \zeta^g ~\bot ~g(\bar{x},\bar{y})\leq 0,~0\le \zeta^u ~\bot ~\bar{u} \geq 0,\qquad\qquad\qquad\qquad\qquad\qquad\qquad\qquad\qquad\qquad\qquad\qquad\qquad\\
\label{WDPKKT-8}&&\nabla_z L(\bar{x},\bar{z}, \bar{u},\bar{v})=0,~h(\bar{x},\bar{y})=0.
\end{eqnarray}
By setting
$\eta^g=\zeta^g,~\eta^h=\zeta^h,~\eta^u=\zeta^u,~\alpha=\widehat{\alpha},~\beta=\widehat{\beta},$
$\xi^g=(\widehat{\alpha},\cdots,\widehat{\alpha})^T\in\mathbb{R}^{p}$, and $\xi^h=-(\widehat{\alpha},\cdots,\widehat{\alpha})^T\in\mathbb{R}^{q}$, \eqref{WDPKKT-1}-\eqref{WDPKKT-8} can be rewritten equivalently as
\begin{eqnarray}
\label{eMDPKKT-1}&&\nabla_x  F(\bar{x},\bar{y})  + \alpha(\nabla_x f(\bar{x},\bar{y})-\nabla_x f(\bar{x},\bar{z}))-
\textstyle\sum\limits_{i=1}^{p}\nabla_x g_i(\bar{x},\bar{z})\bar{u}_i\xi_i^g+\sum\limits_{i=1}^{q}\nabla_x h_i(\bar{x},\bar{z})\bar{v}_i\xi_i^h\nonumber \\
&&+\nabla_{zx}^2 L(\bar{x},\bar{z}, \bar{u},\bar{v}) \beta  +\nabla_x g(\bar{x},\bar{y}) \eta^g+\nabla_x h(\bar{x},\bar{y}) \eta^h=0,  \\
\label{eMDPKKT-2}&&\nabla_y F(\bar{x},\bar{y}) + \alpha\nabla_y f(\bar{x},\bar{y}) + \nabla_y g(\bar{x},\bar{y})\eta^g+ \nabla_y h(\bar{x},\bar{y})\eta^h=0, \\
\label{eMDPKKT-3}&&-\alpha\nabla_z f(\bar{x},\bar{z})+\nabla_{zz}^2 L(\bar{x},\bar{z}, \bar{u},\bar{v}) \beta -\sum\limits_{i=1}^{p}\nabla_z g_i(\bar{x},\bar{z}) u_i\xi_i^g+\sum\limits_{i=1}^{q}\nabla_x h_i(\bar{x},\bar{z}) v_i\xi_i^h =0, \\
\label{eMDPKKT-4}&&- {\rm diag}(g_1(\bar{x},\bar{z}),\cdots,g_p(\bar{x},\bar{z}))\xi^g
+\nabla_z g(\bar{x},\bar{z})^T\beta - \eta^u=0, \\
\label{eMDPKKT-5}&&{\rm diag}(h_1(\bar{x},\bar{z}),\cdots,h_q(\bar{x},\bar{z}))\xi^h+\nabla_z h(\bar{x},\bar{z})^T\beta=0, \\
\label{eMDPKKT-6}&&0\le \eta^g ~\bot ~g(\bar{x},\bar{y})\leq 0, ~0\le \eta^u ~\bot ~\bar{u} \geq 0,\\
\label{eMDPKKT-7}&&\nabla_z L(\bar{x},\bar{z}, \bar{u},\bar{v})=0,~h(\bar{x},\bar{y})=0.
\end{eqnarray}
By \eqref{WDPKKT-6}, we further have
\begin{eqnarray}\label{WDPKKTthm-1}
\widehat{\alpha}(f(\bar{x},\bar{y})-L(\bar{x},\bar{z},\bar{u},\bar{v}))=
\alpha(f(\bar{x},\bar{y})-f(\bar{x},\bar{z}))-\alpha (\bar{u}^Tg(\bar{x},\bar{z})+\bar{v}^Th(\bar{x},\bar{z}))=0.
\end{eqnarray}
Since $(\bar{x},\bar{y},\bar{z},\bar{u},\bar{v})$ is feasible to eMDP, we have $f(\bar{x},\bar{y})-f(\bar{x},\bar{z})\leq0$, $\bar{u}\circ g(\bar{x},\bar{z})\geq0$, and $\bar{v}\circ h(\bar{x},\bar{z})=0$, which implies
\begin{eqnarray}
\label{eMDPKKT-8}&&0\leq\alpha ~\bot ~f(\bar{x},\bar{y})- f(\bar{x},\bar{z})\leq 0, \\
\label{eMDPKKT-9}&&0\le \alpha ~\bot ~\bar{u}^{T} g(\bar{x},\bar{z}) \geq0
\end{eqnarray}
by the nonnegativity of $\alpha$. Moreover, it is easy to verify that \eqref{eMDPKKT-9} is equal to
\begin{eqnarray}
\label{eMDPKKT-10}0\le \xi^g ~\bot ~\bar{u}\circ g(\bar{x},\bar{z}) \geq0.
\end{eqnarray}
Thus, \eqref{eMDPKKT-1}-\eqref{eMDPKKT-7} and \eqref{eMDPKKT-8}-\eqref{eMDPKKT-10}
mean that $(\bar{x}, \bar{y}, \bar{z}, \bar{u},\bar{v})$ is a KKT point of {\rm eMDP}.

(ii) We take the case of eMDP as an example since the other two cases TMDP/eTMDP can be shown in a similar way.
In fact, if $(\bar{x},\bar{y},\bar{z},\bar{u},\bar{v})$ is a KKT point of {\rm MDP}, there is $(\zeta^g, \zeta^h, \zeta^u, \widehat{\alpha}, \widehat{\beta},\widehat{\gamma}) \in \mathbb{R}^{2p+q+2+m}$ such that
\begin{eqnarray}
\label{MDPKKT2-1}&&\nabla_x  F(\bar{x},\bar{y})  + \widehat{\alpha}(\nabla_x f(\bar{x},\bar{y})-\nabla_x f(\bar{x},\bar{z}))-\widehat{\gamma}(\nabla_x g(\bar{x},\bar{z})\bar{u}+\nabla_x h(\bar{x},\bar{z})\bar{v}) \nonumber \\
&&+\nabla_{zx}^2 L(\bar{x},\bar{z}, \bar{u},\bar{v})\widehat{ \beta}  +\nabla_x g(\bar{x},\bar{y}) \zeta^g+\nabla_x h(\bar{x},\bar{y}) \zeta^h=0,  \\
\label{MDPKKT2-2}&&\nabla_y F(\bar{x},\bar{y}) + \widehat{\alpha}\nabla_y f(\bar{x},\bar{y}) + \nabla_y g(\bar{x},\bar{y})\zeta^g+ \nabla_y h(\bar{x},\bar{y})\zeta^h=0, \\
\label{MDPKKT2-3}&&-\widehat{\alpha}\nabla_z f(\bar{x},\bar{z})-\widehat{\gamma}(\nabla_z g(\bar{x},\bar{z})\bar{u}+\nabla_z h(\bar{x},\bar{z})\bar{v})+\nabla_{zz}^2 L(\bar{x},\bar{z}, \bar{u},\bar{v}) \widehat{\beta}=0, \\
\label{MDPKKT2-4}&&-\widehat{\gamma} g(\bar{x},\bar{z})+\nabla_y g(\bar{x},\bar{z})^T\widehat{\beta} - \zeta^u=0, \\
\label{MDPKKT2-5}&&-\widehat{\gamma} h(\bar{x},\bar{z})+\nabla_z h(\bar{x},\bar{z})^T\widehat{\beta}=0, \\
\label{MDPKKT2-6}&&0\leq\widehat{\alpha} ~\bot ~f(\bar{x},\bar{y})- f(\bar{x},\bar{z})\leq 0, \\
\label{MDPKKT2-7}&&0\le \widehat{\gamma} ~\bot ~\bar{u}^{T} g(\bar{x},\bar{z})+\bar{v}^{T} h(\bar{x},\bar{z}) \geq0,\\
\label{MDPKKT2-8}&&0\le \zeta^g ~\bot ~g(\bar{x},\bar{y})\leq 0,~0\le \zeta^u ~\bot ~\bar{u} \geq 0,\\
\label{MDPKKT2-9}&&\nabla_z L(\bar{x},\bar{z}, \bar{u},\bar{v})=0,~h(\bar{x},\bar{y})=0.
\end{eqnarray}
By letting
$\eta^g=\zeta^g,~\eta^h=\zeta^h,~\eta^u=\zeta^u,~\alpha=\widehat{\alpha},~\beta=\widehat{\beta},$
$\xi^g=(\widehat{\gamma},\cdots,\widehat{\gamma})^T\in\mathbb{R}^{p}$, and $\xi^h=-(\widehat{\gamma},\cdots,\widehat{\gamma})^T\in\mathbb{R}^{q}$, \eqref{MDPKKT2-1}-\eqref{MDPKKT2-9} are equivalent to
\begin{eqnarray}
\label{eMDPKKT2-1}&&\nabla_x  F(\bar{x},\bar{y})  + \alpha(\nabla_x f(\bar{x},\bar{y})-\nabla_x f(\bar{x},\bar{z}))-
\textstyle\sum\limits_{i=1}^{p}\nabla_x g_i(\bar{x},\bar{z})\bar{u}_i\xi_i^g+\sum\limits_{i=1}^{q}\nabla_x h_i(\bar{x},\bar{z})\bar{v}_i\xi_i^h\nonumber \\
&&+\nabla_{zx}^2 L(\bar{x},\bar{z}, \bar{u},\bar{v}) \beta  +\nabla_x g(\bar{x},\bar{y}) \eta^g+\nabla_x h(\bar{x},\bar{y}) \eta^h=0,  \\
\label{eMDPKKT2-2}&&\nabla_y F(\bar{x},\bar{y}) + \alpha\nabla_y f(\bar{x},\bar{y}) + \nabla_y g(\bar{x},\bar{y})\eta^g+ \nabla_y h(\bar{x},\bar{y})\eta^h=0, \\
\label{eMDPKKT2-3}&&-\alpha\nabla_z f(\bar{x},\bar{z})+\nabla_{zz}^2 L(\bar{x},\bar{z}, \bar{u},\bar{v}) \beta -\sum\limits_{i=1}^{p}\nabla_z g_i(\bar{x},\bar{z}) u_i\xi_i^g+\sum\limits_{i=1}^{q}\nabla_x h_i(\bar{x},\bar{z}) v_i\xi_i^h =0, \\
\label{eMDPKKT2-4}&&- {\rm diag}(g_1(\bar{x},\bar{z}),\cdots,g_p(\bar{x},\bar{z}))\xi^g
+\nabla_z g(\bar{x},\bar{z})^T\beta - \eta^u=0, \\
\label{eMDPKKT2-5}&&{\rm diag}(h_1(\bar{x},\bar{z}),\cdots,h_q(\bar{x},\bar{z}))\xi^h+\nabla_z h(\bar{x},\bar{z})^T\beta=0, \\
\label{eMDPKKT2-6}&&0\leq\alpha ~\bot ~f(\bar{x},\bar{y})- f(\bar{x},\bar{z})\leq 0, \\
\label{eMDPKKT2-7}&&0\le \eta^g ~\bot ~g(\bar{x},\bar{y})\leq 0, ~0\le \eta^u ~\bot ~\bar{u} \geq 0,\\
\label{eMDPKKT2-8}&&\nabla_z L(\bar{x},\bar{z}, \bar{u},\bar{v})=0,~h(\bar{x},\bar{y})=0.
\end{eqnarray}
By \eqref{MDPKKT2-7}, we have
$0\le \widehat{\gamma} ~\bot ~\bar{u}^{T} g(\bar{x},\bar{z})+\bar{v}^{T} h(\bar{x},\bar{z}) \geq0
$. Since $(\bar{x},\bar{y},\bar{z},\bar{u},\bar{v})$ is feasible to eMDP, we have $h(\bar{x},\bar{z})=0$, which implies
\begin{eqnarray}
\label{eMDPKKT2-9}0\le \xi^g ~\bot ~\bar{u}\circ g(\bar{x},\bar{z}) \geq0.
\end{eqnarray}
\eqref{eMDPKKT2-1}-\eqref{eMDPKKT2-9}
mean that $(\bar{x}, \bar{y}, \bar{z}, \bar{u},\bar{v})$ is a KKT point of {\rm eMDP}.

(iii) This result is almost the same to (ii) and so we omit its proof here.

\vspace{2mm}

\noindent {\bf Proof of Theorem 3.3}
For a feasible point $(\bar{x},\bar{y},\bar{u},\bar{v})$ of MPCC, we define the following index sets:
\begin{eqnarray*}
&&I_{0+}=\{i : g_{i}(\bar{x},\bar{y})=0, ~\bar{u}_i>0\},~~
I_{-0}=\{i : g_{i}(\bar{x},\bar{y})<0, ~\bar{u}_i=0\},~~
I_{00}=\{i : g_{i}(\bar{x},\bar{y})=0, ~\bar{u}_i=0\}.
\end{eqnarray*}	
Recall that $(\bar{x},\bar{y},\bar{u},\bar{v})$ is strongly stationary to {\rm MPCC}
if there is $(\lambda^g,\lambda^h,\lambda^u,\lambda^L)\in\mathbb{R}^{2p+q+m}$
such that
\begin{eqnarray}
&& \nabla_x F(\bar{x},\bar{y}) + \nabla_{yx}^2 L(\bar{x},\bar{y}, \bar{u},\bar{v}) \lambda^L  + \nabla_x g(\bar{x},\bar{y}) \lambda^g+\nabla_x h(\bar{x},\bar{y}) \lambda^h=0,\label{MPCC-1}\\
&&\nabla_y F(\bar{x},\bar{y}) + \nabla_{yy}^2 L(\bar{x},\bar{y}, \bar{u},\bar{v})\lambda^L + \nabla_y g(\bar{x},\bar{y}) \lambda^g+\nabla_y h(\bar{x},\bar{y}) \lambda^h=0, \label{MPCC-2}\\
&&\nabla_y g(\bar{x},\bar{y})^T \lambda^L -\lambda^u= 0, \label{MPCC-3}\\
&&\nabla_y h(\bar{x},\bar{y})^T \lambda^L = 0, \label{MPCC-4}\\
&& \nabla_{y} L(\bar{x},\bar{y}, \bar{u},\bar{v}) =0, ~g(\bar{x},\bar{y})\le 0, ~\bar{u}\ge 0,~h(\bar{x},\bar{y})=0, \label{MPCC-5}\\
&&
\lambda^g_i = 0, ~~i \in I_{-0}, \label{MPCC-6}\\
&&
\lambda^u_i = 0, ~~i \in I_{0+}, \label{MPCC-7}\\
&&\lambda^g_i \geq 0, ~\lambda^u_i \geq 0, ~~i\in I_{00}. \label{MPCC-8}
\end{eqnarray}

Note that the cases of WDP and MDP have been shown in Theorem 4.2 of Li et al. (2023) and Theorem 10 of Li et al. (2024), respectively. Actually, the other cases of eMDP/TWDP/TMDP/eTMDP can be shown in a similar way. For completeness, we give a proof for the case of eMDP as an example.
In fact, if $(\bar{x},\bar{y},\bar{y},\bar{u},\bar{v})$ is a KKT point of eMDP, there exists
$(\eta^g, \eta^h, \eta^u, \alpha, \beta,\xi^g,\xi^h) \in \mathbb{R}^{3p+2q+1+m}$ such that
\begin{eqnarray}
\label{eMDPKKT3-1}&&\nabla_x  F(\bar{x},\bar{y})
+\nabla_{yx}^2 L(\bar{x},\bar{y}, \bar{u},\bar{v}) \beta
+\textstyle\sum\limits_{i=1}^{p}\nabla_x g_i(\bar{x},\bar{y})(\eta_i^g-\bar{u}_i\xi_i^g)+\sum\limits_{i=1}^{q}\nabla_x h_i(\bar{x},\bar{y})(\eta_i^h+\bar{v}_i\xi_i^h)=0,  \qquad\\
\label{eMDPKKT3-2}&&\nabla_y F(\bar{x},\bar{y}) + \alpha\nabla_y f(\bar{x},\bar{y}) + \nabla_y g(\bar{x},\bar{y})\eta^g+ \nabla_y h(\bar{x},\bar{y})\eta^h=0, \\
\label{eMDPKKT3-3}&&-\alpha\nabla_y f(\bar{x},\bar{y})+\nabla_{yy}^2 L(\bar{x},\bar{y}, \bar{u},\bar{v}) \beta -\sum\limits_{i=1}^{p}\nabla_y g_i(\bar{x},\bar{y}) u_i\xi_i^g+\sum\limits_{i=1}^{q}\nabla_y h_i(\bar{x},\bar{y}) v_i\xi_i^h =0, \\
\label{eMDPKKT3-4}&&- {\rm diag}(g_1(\bar{x},\bar{y}),\cdots,g_p(\bar{x},\bar{y}))\xi^g
+\nabla_y g(\bar{x},\bar{y})^T\beta - \eta^u=0, \\
\label{eMDPKKT3-5}&&{\rm diag}(h_1(\bar{x},\bar{y}),\cdots,h_q(\bar{x},\bar{y}))\xi^h+\nabla_y h(\bar{x},\bar{y})^T\beta=0, \\
\label{eMDPKKT3-6}&&0\le \eta^g ~\bot ~g(\bar{x},\bar{y})\leq 0, \\
\label{eMDPKKT3-7}&&0\le \eta^u ~\bot ~\bar{u} \geq 0,\\
\label{eMDPKKT3-8}&&\nabla_z L(\bar{x},\bar{y}, \bar{u},\bar{v})=0,~h(\bar{x},\bar{y})=0.
\end{eqnarray}
Adding \eqref{eMDPKKT3-2} and \eqref{eMDPKKT3-3} together yields
\begin{eqnarray}
\nabla_y F(\bar{x},\bar{y})+\nabla_{yy}^2 L(\bar{x},\bar{y}, \bar{u},\bar{v}) \beta + \textstyle\sum\limits_{i=1}^{p}\nabla_y g_i(\bar{x},\bar{y})(\eta_i^g-\bar{u}_i\xi_i^g)+
\sum\limits_{i=1}^{q}\nabla_y h_i(\bar{x},\bar{y})(\eta_i^h+\bar{v}_i\xi_i^h)=0. \label{eMDPKKT3-9}
\end{eqnarray}
Set
$\lambda^g=\eta^g-\bar{u}\circ\xi^g, ~\lambda^h=\eta^h+\bar{v}\circ\xi^h,~\lambda^u=\eta^u+ \xi^g\circ g(\bar{x},\bar{y})$, and $\lambda^L=\beta$.
We obtain \eqref{MPCC-1}-\eqref{MPCC-5} from \eqref{eMDPKKT3-1} and \eqref{eMDPKKT3-4}-\eqref{eMDPKKT3-9} immediately. Next, we show \eqref{MPCC-6}-\eqref{MPCC-8}.
\begin{itemize}
\item If $i\in I_{-0}$, which means $g_{i}(\bar{x},\bar{y})<0$ and $\bar{u}_i=0$, we have $\eta_i^g=0$ by \eqref{eMDPKKT3-6} and hence $\lambda_i^g=\eta_i^g-\bar{u}_i\xi_i^g=0$. Namely, \eqref{MPCC-6} holds.

\item If $i\in I_{0+}$, which means $g_{i}(\bar{x},\bar{y})=0$ and $\bar{u}_i>0$, we have $\eta_i^u=0$ by \eqref{eMDPKKT3-7} and hence $\lambda_i^u=\eta_i^u+\xi_i^g g_{i}(\bar{x},\bar{y})=0$. Namely, \eqref{MPCC-7} holds.

\item If $i\in I_{00}$, which means $g_{i}(\bar{x},\bar{y})=\bar{u}_i=0$, we have $\lambda_i^g=\eta_i^g\geq 0$ and $\lambda_i^u=\eta_i^u \geq 0$ by \eqref{eMDPKKT3-6} and \eqref{eMDPKKT3-7} respectively. Namely, \eqref{MPCC-8} holds.
\end{itemize}
This completes the proof.

\setcounter{equation}{0}
\renewcommand\theequation{C.\arabic{equation}}
\setcounter{thm}{0}
\renewcommand\thethm{C.\arabic{thm}}
\setcounter{defi}{0}
\renewcommand\thedefi{C.\arabic{defi}}
\setcounter{table}{0}
\renewcommand\thetable{C.\arabic{table}}
\setcounter{ex}{0}
\renewcommand\theex{C.\arabic{ex}}

\section*{Appendix C: Proof of Theorem 4.2}
We only show the case of $\mathrm{TWDP}(t)$ because the cases of $\mathrm{TMDP}(t)$/$\mathrm{eTMDP}(t)$ can be shown in a similar way. For simplicity, we take the upper constraint $(x,y)\in \Omega$ away from $\mathrm{TWDP}(t)$. In addition, we denote
\begin{eqnarray*}
&&\bar{w}=(\bar{x},\bar{y},\bar{z},\bar{u},\bar{v}), ~~~w^k=(x^k,y^k,z^k,u^k,v^k), ~~~w=(x,y,z,u,v),\\
&&F_0(w)=F(x,y),~~~G_0(w)=f(x,y)-f(x, z)-u^Tg(x,z),\\
&&G_i(w)=g_i(x,y)~(1\le i\le p),~~~G_i(w)=-u_{i-p} ~(p+1\le i\le 2p),~~~H_i(w)=h_i(x,y) ~(1\le i\le q),\\
&& H_i(w)=h_{i-q}(x,z) ~(q+1\le i\le 2q), ~~H_i(w)=\nabla_{z_{i-2q}}L(x,z,u,v) ~(2q+1\le i\le 2q+m).
\end{eqnarray*}
Then, TWDP and $\mathrm{TWDP}(t_k)$ become
\begin{eqnarray}\label{thmRelaxation-1}
\min &&F_0(w) \nonumber \\
\mbox{\rm s.t.} && G_i(w) \leq 0 ~(i=0,1,\cdots,2p),\qquad \qquad\qquad\qquad\\
&&H_i(w) =0 ~(i=1,\cdots,2q+m) \nonumber
\end{eqnarray}
and
\begin{eqnarray}\label{thmRelaxation-2}
\min &&F_0(w) \nonumber \\
\mbox{\rm s.t.} && G_0(w) \leq t_k, \\
&& G_i(w) \leq 0 ~(i=1,\cdots,2p),\ \ H_i(w) =0 ~(i=1,\cdots,2q+m). \nonumber
\end{eqnarray}
Obviously, $\bar{w}$ is a feasible point of TWDP and, without loss of generality, we may assume that the whole sequence $\{w^k\}$ converges to $\bar{w}$.
Since $w^k$ is a KKT point of \eqref{thmRelaxation-2}, it follows from Lemma A.1 in Steffensen and Ulbrich (2010) that there exists $(\lambda_0^k,\lambda^k,\mu^k) \in \mathbb{R}\times\mathbb{R}^{2p+1}\times\mathbb{R}^{q+m}$ such that
\begin{eqnarray}
&&\nabla F_0(w^k)+\lambda_0^k\nabla G_0(w^k)+\nabla G(w^k)\lambda^k+\nabla H(w^k)\mu^k=0,
\nonumber\\
&&\lambda_0^k(G_0(w^k)-t_k)=0, ~~\lambda_0^k\geq0,
\label{thmRelaxation-3}\\
&&G(w^k)^T\lambda^k=0, ~~\lambda^k\geq0,\nonumber
\end{eqnarray}
and the gradients
\begin{eqnarray}
\{\nabla G_i(w^k): 
\lambda_i^k>0, 0\le i\le 2p\}\cup\{\nabla H_i(w^k): 
\mu_i^k\not=0\} \label{independent}
\end{eqnarray}
are linearly independent. To show that $\bar{w}$ is a KKT point of \eqref{thmRelaxation-1}, it is sufficient to show the boundedness of the sequence $\{\lambda_0^k,\lambda^k,\mu^k\}$.

Assume by contradiction that $\{\lambda_0^k,\lambda^k,\mu^k\}$ is unbounded. Taking a subsequence if necessary, we may assume
\begin{eqnarray}\label{thmRelaxation-4}
\lim_{k\to\infty}\frac{(\lambda_0^k,\lambda^k,\mu^k)}
{\|(\lambda_0^k,\lambda^k,\mu^k)\|} = (\lambda_0^*,\lambda^*,\mu^*). 
\end{eqnarray}
In particular, for every $k$ sufficiently large, we have
\begin{eqnarray}
\lambda_i^*>0 \ \Rightarrow \ \lambda_i^k>0, \quad \mu_i^*\not=0 \ \Rightarrow \ \mu_i^k\not= 0. \label{supp}
\end{eqnarray}
Dividing by $\|(\lambda_0^k,\lambda^k,\mu^k)\|$ in \eqref{thmRelaxation-3} and taking a limit, we have
\begin{eqnarray}\label{thmRelaxation-5}
&&\lambda_0^*\nabla G_0(\bar{w})+\sum_{\lambda_i^*>0}\lambda_i^*\nabla G_i(\bar{w})+\sum_{\mu_i^*\not=0}\mu_i^*\nabla H_i(\bar{w})=0, ~~\lambda_0^*\geq0.
\end{eqnarray}
Since $\|(\lambda_0^*,\lambda^*,\mu^*)\|=1$ by \eqref{thmRelaxation-4}, \eqref{thmRelaxation-5} implies that the gradients
\begin{eqnarray*}
\{\nabla G_i(\bar{w}): \lambda_i^*>0, 0\le i\le 2p\}\cup\{\nabla H_i(\bar{w}): \mu_i^*\not=0\}
\end{eqnarray*}
are positive-linearly dependent. Since problem \eqref{thmRelaxation-1} satisfies the CPLD in Qi and Wei (2000) at $\bar{w}$, it is not difficult to see that the gradients
$\{\nabla G_i(w): \lambda_i^*>0, 0\le i\le 2p\}\cup\{\nabla H_i(w): \mu_i^*\not=0\}$
are linearly dependent in a neighborhood. This, together with \eqref{supp}, contradicts the linear independence of the vectors in \eqref{independent}.
Consequently, $\{\lambda_0^k,\lambda^k,\mu^k\}$ is bounded.
This completes the proof.

\end{appendices}


\begin{thebibliography}{99}


\bibitem{Ben1990computational}
Ben-Ayed O, Blair C E.
Computational difficulties of bilevel linear programming. Operations Research, 1990, 38(3): 556--560.

\bibitem{Byeon}
Byeon G, Van Hentenryck P. Benders subproblem decomposition for bilevel problems with convex
follower. INFORMS Journal on Computing, 2022, 34(3): 1749--1767.

\bibitem{Caprara}
Caprara A, Carvalho M, Lodi A, Woeginger G J. Bilevel knapsack with interdiction constraints.
INFORMS Journal on Computing, 2016, 28(2): 319--333.


\bibitem{Colson2007overview}
Colson B, Marcotte P, Savard G.
An overview of bilevel optimization. Annals of Operations Research, 2007, 153: 235--256.


\bibitem{Dempe2020bilevel}
Dempe S, Zemkoho A.
Bilevel Optimization: Advances and Next Challenges, vol. 161. Springer Optimization and its Applications, Berlin, 2020.


\bibitem{Goyal2023decision}
Goyal A, Zhang Y, He C.
Decision rule approaches for pessimistic bilevel linear programs under moment ambiguity with facility location applications. INFORMS Journal on Computing, 2023, 35(6): 1342--1360.


\bibitem{LinEquiation}
Guo L, Lin G H, Ye J J. Solving mathematical programs with equilibrium constraints. Journal of Optimization Theory and Applications, 2015, 166(1), 234--256.

\bibitem{Hong2023two}
Hong M, Wai H T, Wang Z, Yang Z.
A two-timescale stochastic algorithm framework for bilevel optimization: Complexity analysis and application to actor-critic. SIAM Journal on Optimization, 2023, 33(1): 147--180.


\bibitem{Izmailov2012semismooth}
Izmailov A F, Pogosyan A L, Solodov M V.
Semismooth Newton method for the lifted reformulation of mathematical programs with complementarity constraints. Computational Optimization and Applications, 2012, 51(1): 199--221.


\bibitem{Khorramfar2022managing}
Khorramfar R, \"{o}zalt{\i}n O Y, Kempf K G, Uzsoy R. Managing product transitions: A bilevel programming approach. INFORMS Journal on Computing, 2022, 34(5): 2828--2844.


\bibitem{Kleinert2021computing}
Kleinert T, Schmidt M.
Computing feasible points of bilevel problems with a penalty alternating direction method. INFORMS Journal on Computing, 2021, 33(1): 198--215.


\bibitem{Leyffer2006interior}
Leyffer S, L\'{o}pez-Calva G, Nocedal J.
Interior methods for mathematical programs with complementarity constraints. SIAM Journal on Optimization, 2006, 17(1): 52--77.


\bibitem{Li2023novel}
Li Y W, Lin G H, Zhang J, Zhu X.
A novel approach for bilevel programs based on Wolfe duality. 2023, https://doi.org/10.48550/arXiv.2302.06838.


\bibitem{Li-PJO}
Li Y W, Lin G H, Zhu X.
Notes on lower-level duality approach for bilevel programs. Pacific Journal of Optimization, 2024, 20(3): 475--488.


\bibitem{Li2024solving}
Li Y W, Lin G H, Zhu X.
Solving bilevel programs based on lower-level Mond-Weir duality. INFORMS Journal on Computing, 2024, 36(5): 1225--1241.


\bibitem{Lin2009solving}
Lin G H, Chen X, Fukushima M.
Solving stochastic mathematical programs with equilibrium constraints via approximation and smoothing implicit programming with penalization.  Mathematical Programming, 2009, 116(1): 343--368.


\bibitem{Lin2006hybrid}
Lin G H, Fukushima M.
Hybrid approach with active set identification for mathematical programs with complementarity constraints.
Journal of Optimization Theory and Applications, 2006, 128(1): 1--28.


\bibitem{Lin2014simple}
Lin G H, Xu M, Ye J J.
On solving simple bilevel programs with a nonconvex lower level program.
Mathematical Programming, 2014, 144(1): 277--305.


\bibitem{Liu2023hierarchical}
Liu R, Liu X, Zeng S, Zhang J, Zhang Y.
Hierarchical optimization-derived learning. IEEE Transactions on Pattern Analysis and Machine Intelligence, 2023, 45(12): 14693--14708.



\bibitem{Lu2024first}
Lu Z, Mei S.
First-order penalty methods for bilevel optimization. SIAM Journal on Optimization, 2024, 34(2): 1937--1969.


\bibitem{Luo1996mathematical}
Luo Z Q, Pang J S, Ralph D.
Mathematical Programs with Equilibrium Constraints.
Cambridge University Press, Cambridge (1996).


\bibitem{Mehlitz2021sufficient}
Mehlitz P, Zemkoho A B.
Sufficient optimality conditions in bilevel programming. Mathematics of Operations Research, 2021, 46(4): 1573--1598.



\bibitem{Mond1981generalized}
Mond B, Weir T.
Generalized concavity and duality, in ``Generalized Concavity in
Optimization and Economics'' (Schaible S and Ziemba W T, Eds.), 263--279,
Academic Press, San Diego (1981).


\bibitem{Qi2000constant}
Qi L, Wei Z.
On the constant positive linear dependence condition and its application to SQP methods. SIAM Journal on Optimization, 2000, 10(4): 963--981.


\bibitem{Sabach2017first}
Sabach S, Shtern S.
A first order method for solving convex bilevel optimization problems. SIAM Journal on Optimization, 2017, 27(2): 640--660.


\bibitem{Scholtes2001}
Scholtes S.
Convergence properties of a regularization scheme for mathematical programs with complementarity constraints. SIAM Journal on Optimization, 2001, 11(4): 918--936.


\bibitem{Steffensen2010new}
Steffensen S, Ulbrich M.
A new relaxation scheme for mathematical programs with equilibrium constraints. SIAM Journal on Optimization, 2010, 20(5): 2504--2539.


\bibitem{Wolfe1961duality}
Wolfe P.
A duality theorem for nonlinear programming.  Quarterly of Applied Mathematics, 1961, 19(3): 239--244.


\bibitem{Xu2014smoothing}
Xu M, Ye J J.
A smoothing augmented Lagrangian method for solving simple bilevel programs. Computational Optimization and Applications, 2014, 59(1--2): 353--377.


\bibitem{Ye2023difference}
Ye J J, Yuan X, Zeng S, Zhang J.
Difference of convex algorithms for bilevel programs with applications in hyperparameter selection. Mathematical Programming, 2023, 198(2): 1583--1616.


\bibitem{Ye2010new}
Ye J J, Zhu D.
New necessary optimality conditions for bilevel programs by combining the MPEC and value function approaches. SIAM Journal on Optimization, 2010, 20(4): 1885--1905.


\bibitem{Zhang2024introduction}
Zhang Y, Khanduri P, Tsaknakis I, Yao Y, Hong M, Liu S.
An introduction to bilevel optimization: Foundations and applications in signal processing and machine learning. IEEE Signal Processing Magazine, 2024, 41(1): 38--59.


\bibitem{Zeng2020practical}
Zeng B. A practical scheme to compute the pessimistic bilevel optimization problem. INFORMS Journal on Computing, 2020, 32(4): 1128--1142.

\end{thebibliography}
\end{document}